\documentclass{amsart}

\usepackage{amssymb}
\usepackage{latexsym}
\usepackage{amsmath}
\usepackage{enumerate}
\usepackage{amsthm}
\usepackage{wasysym}
\usepackage{stmaryrd}
\usepackage{mathrsfs}
\usepackage{pifont}
\usepackage{tikz}
\usepackage{bbm}

\newtheorem{theorem}{Theorem}[section]
\newtheorem{define}[theorem]{Definition}

\newtheorem{exa}[theorem]{Example}

\newtheorem{exerc}[theorem]{Exercise}

\newtheorem{conj}[theorem]{Conjecture}

\newtheorem{ques}[theorem]{Open Question}
\newenvironment{question}{\begin{ques} \rm}{\qee\end{ques}}
\newtheorem{lem}[theorem]{Lemma}
\newenvironment{lemma}{\begin{lem} \it}{\end{lem}}
\newtheorem{cor}[theorem]{Corollary}
\newenvironment{corollary}{\begin{cor} \it}{\end{cor}}
\newtheorem{rem}[theorem]{Remark}
\newenvironment{remark}{\begin{rem} \rm}{\qee\end{rem}}
\newtheorem{obs}[theorem]{Observation}
\newenvironment{observation}{\begin{obs} \it}{\end{obs}}
\newtheorem{athm}[theorem]{Announced Theorem}
\newenvironment{antheorem}{\begin{athm} \it}{\end{athm}}

\definecolor{uuxgreen}{cmyk}{1,0,0.75,0}
\definecolor{uuxred}{cmyk}{0.2,1,0.9,0.1}
\definecolor{uuyblue}  {cmyk}{0.9,0.55,0,0}

\newcommand{\qee} {\hspace*{2mm}\hfill \ding{109}}

\renewcommand{\iff}{\leftrightarrow}
\renewcommand{\phi}{\varphi}

\DeclareMathOperator{\possible}{\text{\tikz[scale=.6ex/1cm,baseline=-.6ex,rotate=45,line width=.1ex]{
                            \draw (-1,-1) rectangle (1,1);}}}
\DeclareMathOperator{\necessary}{\text{\tikz[scale=.6ex/1cm,baseline=-.6ex,line width=.1ex]{
                            \draw (-1,-1) rectangle (1,1);}}}

 \DeclareMathOperator{\gnecessary}{\text{\tikz[scale=.6ex/1cm,baseline=-.6ex,line width=.1ex]{
                            \draw[gray, fill = gray, fill opacity = .90] (-1,-1) rectangle (1,1);}}}
                            
\DeclareMathOperator{\dotnecessary}{\text{\tikz[scale=.6ex/1cm,baseline=-.6ex,line width=.1ex]{
                            \draw (-1,-1) rectangle (1,1);  \draw[fill=black] (0,0) circle (.25);}}}
                     
\newcommand{\qedright}{\belowdisplayskip=-12pt}
\newcommand{\grullet}{\text{\textcolor{gray}{$\bullet$}}}
\newcommand{\subun}[1]{[#1]}
\newcommand{\sbra}[1]{\textup(#1\textup)}
\newcommand{\grafo}{{\mathbbm g}}
\newcommand{\sugr}[2]{{#1\downarrow #2}}
\newcommand{\snip}[2]{(#1\mathbin{\lightning} #2)}
\newcommand{\bisim}{\simeq}
\newcommand{\iso}{\cong}

\newcommand{\eva}[1]{\llbracket #1 \rrbracket}
 \newcommand{\tupel}[1]{{\langle #1 \rangle}}
\newcommand{\verz}[1]{\{ #1 \}}
\newcommand{\To}{\Rightarrow}
\newcommand{\Iff}{\Leftrightarrow}

\newcommand{\gnum}[1]{{\ulcorner #1 \urcorner}}
\newcommand{\apr}{{\vartriangle}}

\newcommand{\opr}{\necessary}
\newcommand{\dotbox}{\dotnecessary}

\newcommand{\oco}{\possible}

\newcommand{\graysq}{\gnecessary}
\newcommand{\glc}{\ensuremath{{\sf GL}^\circ}}

\title[Cyclic Henkin Logic]{Cyclic Henkin Logic}



\author{Albert Visser}
 \address{Philosophy, Faculty of Humanities,
                Utrecht University,
               Janskerkhof 13,
                3512BL~~Utrecht, The Netherlands}
\email{a.visser@uu.nl}
\date{\today}

\thanks{I am grateful to Johan van Benthem and to the participants
of the Implicit Seminar Warsaw, the Proof Theory Virtual Seminar, and the Logic Online Seminar Moscow
for their comments on and questions about various talks I gave on the subject of this paper or related matters.
I thank Tadeusz Litak with whom I am writing a parallel paper. Without our shared project
the present paper would not be possible.}

\begin{document}

\begin{abstract}
In this paper, we study Cyclic Henkin Logic {\sf CHL}, a  logic that can be described
as provability logic without the third L\"ob condition, to wit, that provable implies provably provable (aka principle 4).
The logic {\sf CHL} does have full modalised fixed points. We implement these fixed points using
cyclic syntax, so that we can work just with the usual repertoire of connectives. 

The main part of the paper is devoted to developing the logic on cyclic syntax. 
Many theorems, like the multiple fixed point theorem, become matter of course in this context.
We submit that the use of cyclic syntax 
is of interest even for the study of classical L\"ob's Logic. We show that a version of the
de Jongh-Sambin algorithm can be seen as one half of a synonymy between 
${\sf GL}^\circ$, i.e.\ {\sf CHL} plus the third L\"ob Condition, and
ordinary L\"ob's Logic {\sf GL}. Our development illustrates  that an appropriate computation scheme for the algorithm 
is guard recursion.

We show how arithmetical interpretations work for the cyclic syntax. In an appendix, we 
give some further information about the arithmetical side of the equation.
\end{abstract}

\maketitle

\section{Introduction}
In the present paper, we study Cyclic Henkin Logic {\sf CHL}. Our original interest in this logic was triggered by
 the development, in our paper \cite{viss:look19}, of a class of provability predicates  for which
 we do have L\"ob's Rule but for which the third L\"ob condition, to wit, provable implies provably
 provable (aka principle 4), fails. It turns out that, in this context, the  
  de Jongh-Sambin-Bernardi Theorem about the uniqueness of modalised fixed points does hold. 
  Thus, a substantial amount of familiar reasoning from L\"ob's Logic {\sf GL} is still  present. On the other hand,
  in the absence of the third L\"ob condition, the de Jongh-Sambin Theorem about the
  explicit definability of fixed points fails. In fact, the most salient fixed point of them all, the G\"odel fixed point,
   does not have an explicit definition.
  
  So, what is the appropriate form for an appropriately weakened version of L\"ob's Logic?
Our intuition is that the logic in question should have full modalised fixed points. Given that
  design choice, in the light of non-explicit definability, there seem to be essentially three roads to
  follow: (i) we try to extend the repertoire of ordinary modal operators so that we do recover
  explicit definability for the enriched repertoire; (ii) we add a variable-binding fixed point operator
  $\digamma p.\phi$ that is only allowed when $\phi$ is modalised in $p$; (iii) we work with cyclic syntax.
 We have not looked seriously at option (i), but it very well possible that it is not feasible in a reasonable way.
 Option (ii) will be worked out in a forthcoming paper by Tadeusz Litak and the author. In the present paper
 we study option (iii).
 
 We think that the cyclic syntax has a wider interest than just the study of a weaker version of
 L\"ob's Logic. It also gives us a new way to look at L\"ob's Logic itself. We will show, for example,
 that a version of the de Jongh-Sambin algorithm to compute explicit fixed points can be viewed
 as the specification of one half of a synonymy between ${\sf GL}^\circ$, i.e.\ {\sf CHL} plus the third L\"ob condition, and
 {\sf GL} in its original formulation. The specification of the algorithm shows that guard recursion is an
 appropriate way of thinking about this algorithm.
 
 The logic {\sf CHL} is synonymous with a corresponding theory, the Henkin Calculus or {\sf HC}, which is
 formulated in a language with a variable-binding fixed point operator. This will be proved in a forth-coming paper
 by Tadeusz Litak and the author. In that same paper, we will prove that {\sf HC}, and, thus, also {\sf CHL}, is synonymous
 with the well-founded part of the $\mu$-calculus, i.e., the $\mu$-calculus plus the minimal Henkin sentence
 $\mu p.\opr p$. This last result is based on ideas from Johan van Benthem's paper  \cite{bent:moda06}, which were
 extended in a  paper by the author \cite{viss:lobs05}. Thus, {\sf CHL} can be viewed as a treatment of the well-founded
 part of the $\mu$-calculus on a  quite different syntax.
 
  An obvious further step in the project of provability logic on cyclic syntax is to connect circular syntax with circular proofs as studied in, e.g.,
 \cite{sham:circ14} and \cite{sham:nonw20}. We have not explored this attractive possibility yet.

\subsection{Plan of the Paper} 
 A substantial part of the paper is devoted to carefully developing the system.
 We develop the syntax and introduce the appropriate principles of definition and proof concerning the syntax
  in Section~\ref{syntaxsection}.
 Then, we introduce {\sf CHL} and work out the basic facts about the logic in Section~\ref{cyhelo}.
 In Section~\ref{leulosection}, we prove the synonymy between ${\sf GL}^\circ$ and {\sf GL}.
 We study some further inter-theory relations in Section~\ref{intersection}. 
 Section~\ref{arithintsection} contains our development of arithmetical interpretations. 
 Finally, Appendix~\ref{arima} provides a somewhat closer look at the arithmetical side of the equation.
 
  \subsection{Prerequisites}
Some knowledge of basic provability logic is helpful. The classical textbooks \cite{bool:logi93} or
\cite{smor:self85} are quite sufficient. However, there are many other good expositions available nowadays,
like \cite{lind:prov96}, \cite{japa:logi98}, \cite{svej:prov00}, \cite{arte:prov04}. 

It would also be good if
the reader has at least seen the notion of bisimulation. Any modern introduction to modal logic or to
computer science will explain this notion. The books \cite{barw:liar87} and \cite{acze:nonw88} also
introduce the notion and also contain some material closely connected to the present paper.

The appendix contains some more advanced material on arithmetic, but the rest of the paper is
independent of this. 

\section{Syntax}\label{syntaxsection}
In this Section we provide the basics of our cyclic syntax. A major inspiration of our approach is the treatment of
the paradoxes in the book \cite{barw:liar87}.

In Subsection~\ref{graphs}, we specify what kind of graphs
we employ in our development. Then, in Subsection~\ref{formulas}, we introduce the relevant notion of formula and
develop some basic proof methods and definition methods. Also, we prove a number of basic facts.

\subsection{Graphs}\label{graphs}
Let a non-empty set of labels $L$ be fixed. We have an arity function ${\sf ar}:L \to \omega$.

We need the following notion of graph: a directed pointed labeled graph with ordered successors.
A graph $\mathcal G$ is given as a quadrupel $\tupel{V,r,S,\lambda}$.
 Here:
\begin{itemize}
\item
 $V$ is the set of vertices or nodes. In our paper this set will always be finite.
\item
$r\in V$ is the point or root.
\item
$\lambda: V \to L$ is the labeling function. A vertex $a$ with label $\mathfrak a$ will  be called \emph{an
occurrence of $\mathfrak a$}.
\item
 $S:V\to V^\ast$, where $V^\ast$ is the set of (finite) sequences
of elements in $V$ (including the empty sequence). 
We demand that ${\sf length}(Sa) = {\sf ar}(\lambda a)$ (in $\mathcal G$).

We write $S_i(a)$ for $(S(a))_i$, where $i < {\sf length}(Sa)$.
\item
Let $\hat S$ be the relation given by $a\hat S b $ iff there is an $i< {\sf length}(Sa)$ such that $b= S_ia$.
A \emph{path} in the a graph is a sequence of vertices $a_0\hat S a_1\hat S \dots \hat S a_{k}$. 
We demand that every vertex can be reached
via a finite path from $r$.
\end{itemize}

For many purposes the notion of path as defined here is sufficient. However it is also good to have the
notion of \emph{directive path}. To motivate this consider a graph with nodes $a$ and $b$ and suppose
$Sa = \tupel{b,b}$ and $Sb = \tupel{a}$. We can have a path $aba$ that takes the left turn and one 
 that takes the right turn. To distinguish such possibilities, we define \emph{a directive path} as a sequence
$a_0i_0\dots a_{k-1}i_{k-1}a_k$, where $i_j < {\sf ar}(\lambda a_j)$ and $a_{j+1}=S_{i_j}a_j$.

\begin{remark}
We opted for the present format for graphs since this is in accordance with the representation as co-algebra.
Of course other formats are possible.
\end{remark}

A \emph{bisimulation} between graphs $\mathcal G$ and $\mathcal G'$ is a relation $R$ between
$V$ and $V'$, such that:
\begin{enumerate}[i.]
\item
 If $aRa'$, then $\lambda a=\lambda' a'$.
 \item
If $aRa'$, then $S a$ and $S' a'$ are sequences with the same length $\ell$ and we have $S_i a \mathrel R S'_ia'$, for all $i< \ell$.
\end{enumerate}

\noindent
Two graphs are \emph{bisimilar} if there is a bisimulation between them that relates their roots. 
We note that is follows that the bisimulation is total and surjective.
We write $\bisim$ for bisimilarity.
We remind the reader of the well-known fact that bisimulations are closed under unions.
Thus, there is a maximal bisimulation between any two graphs. Of course, this does not
need to relate the roots.

An isomorphism between graphs is a bijective bisimulation that relates the roots.
We write $\iso$ for isomorphism. We will think of our graphs modulo isomorphism. 
Philosophically, we want to think about the graphs modulo bisimulation, however, it is technically
convenient to have the more `concrete' representations modulo isomorphism available.

\begin{remark}
If we would allow infinite graphs, we could define the canonical unraveling of $\mathcal G$ as the graph
$\mathcal G'$ with as domain the directive paths in $\mathcal G$. The new successor and labeling functions
are as expected. One can then show that two graphs are bisimilar iff their canonical unravelings are isomorphic.
\end{remark}

Rooted graphs are our default. On occasion we will also consider \emph{unrooted graphs}.
Of course these are just graphs minus the root. We will allow such graphs to contain disconnected
parts. The definition of bisimulation remains the same without the condition for the roots.
In stead we demand the relation to be total and surjective.

A \emph{cycle} in a graph is a set of vertices $C$ such that we can arrange the elements of $C$
in a path $a_0\hat S a_1 \hat S \dots \hat Sa_{k-1}\hat S a_0$. 
We demand that the $a_j$ are pairwise distinct.
Note that in our definition a cycle has no designated starting point.
A vertex is a \emph{cycle vertex} if it is on a cycle.

We will write $\mathsf c(\mathcal G)$ for the number of cycles in $\mathcal G$.

\begin{remark}
There is also 
the notion of \emph{directive cycle}. A directive cycle is a function $\gamma$ from
a set of vertices $C$ to numbers, such that $\gamma(a)<{\sf ar}(\lambda a)$.
We demand that we can arrange the elements of $C$ in a directive path
$a_0i_0a_1 \dots a_{k_1}i_{k-1}a_0$, where $i_j = \gamma(a_j)$.  
Here the $a_j$ are pairwise distinct.
\end{remark}

We define a number of operations on graphs.
\begin{itemize}
\item
$\grafo(a,\mathcal G_0,\dots, G_{k-1})$ is the result of taking the disjoint sum of the 
$\mathcal G_0,\dots,\mathcal G_{k-1}$
and adding a fresh root $r$ with label $a$ to this sum. We allow that $k=0$ here. 
\item
$\sugr{\mathcal G}{a} =:\mathcal G'$ is subgraph of $\mathcal G$ generated by $a$. It is defined as follows. 
\begin{itemize}
\item
$V'$ is the set of vertices that can be reached via a (possibly empty) path from $a$.
\item
$r' := a$.
\item
$\lambda'$ is $\lambda$ restricted to $V'$.
\item
$S'$ is $S$ restricted to $V'$.
\end{itemize}
\item 
Suppose the length of $Sr$ is $n$ and $i<n$. Then ${\sf su}_i(\mathcal G) := \sugr{\mathcal G}{S_ir}$. 
\item
${\sf min}(\mathcal G)$ is the result of dividing out the maximal auto-bisimulation of $\mathcal G$. 
\end{itemize}

We have the following obvious lemmas.

\begin{lemma}
Suppose $\mathcal G_0\bisim \mathcal G_1$ via $R$ and $a_0Ra_1$.
Then $\sugr{\mathcal G_0}{a_0}\bisim\sugr{\mathcal G_1}{a_1}$ via
the restriction of $R$ to the nodes of $\sugr{\mathcal G_0}{a_0}$ and
$\sugr{\mathcal G_1}{a_1}$.
\end{lemma}

\begin{lemma}
${\sf min}(\mathcal G)$ is bisimulation minimal: all bisimulations on it are subsets of the identity relation on $V$. Moreover,
$\mathcal G\bisim\mathcal G'$  iff ${\sf min}(\mathcal G)\iso {\sf min}(\mathcal G')$.
\end{lemma}

Consider a graph $\mathcal G$. Let $W\subseteq V$.  We say that $W$ is a \emph{guard} (for $\mathcal G$) if every 
 cycle contains an element of $W$.

Thus, we have \emph{guard-induction} and \emph{guard-recursion} in a guarded graph: 
\begin{lemma}\label{barsmurf}
Suppose $W$ is a guard of $\mathcal G$. We have:
\begin{enumerate}[i.]
\item
Suppose we have a property $P$ of vertices such that all $g\in W$ have $P$. Suppose further that
if $a\not\in W$ and all $a$'s successors have $P$, then $a$ has $P$. Then all vertices have $P$.
\item
Suppose for every label $\mathfrak a$, we have a function $G_{\mathfrak a}:D^{{\sf ar}(\mathfrak a)} \to D$ and suppose 
$F:W \to D$.
Then there is a unique function $H:V \to D$ such that $H(a) = F(a)$ if $a\in W$ and
$H(a) = G_{\mathfrak a}(H(S_0a),\dots,H(S_{k-1}a))$ if $a\not\in W$ and $\lambda a = \mathfrak a$ and $k = {\sf ar}(\mathfrak a)$. 
\end{enumerate}
\end{lemma}

\begin{proof}
Ad (i): Consider any $a$ in $\mathcal G$. Let $B$ be the union of $W$ with the set of leaves of $\mathcal G$.
Consider any vertex $a$. Clearly in $\sugr{\mathcal G}{a}$, the intersection $C$ of $B$ with the vertices in
$\sugr{\mathcal G}{a}$ forms a bar, i.e. every indefinitely prolonged path must eventually pass through an element of $C$.
So by bar-induction, the root $a$ has property $P$. The proof of (ii) is similar using bar recursion.
\end{proof}

Suppose $W$ is a guard for $\mathcal G$ and $W'$ is a guard for $\mathcal G'$.
A bisimulation $R$ between $\mathcal G$ and $\mathcal G'$ is a $W,W'$-bisimulation if 
whenever $aRa'$, then $a\in W$ iff $a'\in W'$.

\begin{lemma}
Suppose for every label $\mathfrak a$, we have a function $G_{\mathfrak a}:D^{{\sf ar}(\mathfrak a)} \to D$.
Let $G$ be a gard for $\mathcal G$ and let $G'$ be a guard for $\mathcal G'$. Suppose $R$ is a $G,G'$-bisimulation
between $\mathcal G$ and $\mathcal G'$. ($R$ does not need to be root-preserving.)
Consider $F:G \to D$ and $F':G'\to D$ such that whenever $a\in W$ and $aRa'$, we have $F(a)=F'(a')$.
Let $H$ and $H'$ be the functions guaranteed to exist by Lemma~\ref{barsmurf}.
Then, for all $a\in V$ and $a'\in V'$, if $aRa'$, then $H(a)=H(a')$.
\end{lemma}

\begin{proof}
Let $P(a)$ be the property: for all $a'$ if $aRa'$, then $H(a) = H(a')$. Clearly we have $P$ on $W$.
Moreover, it is easy to see that $P$ is preserved from the successors of $a$ to $a$.
\end{proof}

\subsection{Formulas}\label{formulas}
We define our formulas. The totality of formulas will constitute our full language $\mathbb L^\circ$.

 The set of labels $L$ for formulas is given by $\bot$, $\top$, $\neg$, $\opr$,
$\wedge$, $\vee$, $\to$, $p_0$, $p_1$, \dots,
where ${\sf ar}(\top)={\sf ar}(\bot) = {\sf ar}(p_i) =0$, ${\sf ar}(\neg)={\sf ar}(\opr) =1$, ${\sf ar}(\wedge)={\sf ar}(\vee)={\sf ar}(\to)=2$. 
As usual we use also $p,q,r,r',\dots$ for propositional variables. 
A \emph{formula} $\phi$ is a  graph for which ${\sf bo}(\phi)$, the set of $\opr$-occurrences, is a guard.
We note that the set of $\opr$-occurrences on a cycle ${\sf bo}^\circ(\phi)$
also forms a guard.

If the label of the root is not a variable, we call it \emph{the main connective}.
We say that \emph{$\phi$ is modalised in $p$} if every path from the root to an occurrence of $p$ contains a $\opr$-occurrence.

\begin{remark}
We can very well model our formulas in the hereditarily finite sets of non-well-founded 
set theory {\sf AFA} with the labels as ur-elements.
Only the guardedness condition is perhaps somewhat unnatural in this context. 
\end{remark}

We define some operations on formulas:
\begin{itemize}
\item
We identify $\top$ with $\grafo(\top)$ and, similarly, for $\bot$ and the $p_i$.\footnote{The convenient confusion between
labels and operations has its limits should be treated with some care. We will warn the reader when to tread carefully.}
\item
$\neg\,\phi := \grafo(\neg,\phi)$ and, similarly, for $\opr$.
\item
$(\phi \wedge \psi) :=  \grafo(\wedge,\phi,\psi)$ and, similarly, for $\vee$ and $\to$.
\item
$\digamma p.\phi$ is the result of identifying the root with all vertices labeled $p$, where one keeps the
label of the root. This is only allowed when $\phi$ is modalised in $p$, since otherwise the resulting graph will
not be a formula.
It is easily seen that in the resulting graph all cycles contain a $\opr$-occurrence, as desired.
\item
$\phi[p_0:\psi_0,\dots,p_{k-1}:\psi_{k-1}]$ is the result of the following operation. First, we form the disjoint union
of  $\phi$ and the $\psi_i$. Then, simultaneously, we identify the vertices labeled $p_i$ in (the disjoint copy of) $\phi$ with the root of
(the disjoint copy of) $\psi_i$,
where we keep the label of the root of $\psi_i$ and discard the label $p_i$.
\end{itemize}

\begin{remark}
Suppose $\phi$ is modalised in $p$ and $q$. Then $\digamma p.\digamma q.\phi \iso \digamma p. \phi [q:p]$.
Thus, we see that one of the costs of the graph approach is that this principle is built in. Of course, one may also
consider it as a bonus.
\end{remark}

The operations are safe for bisimulation:  

\begin{lemma}
Bisimilarity between formulas is a congruence relation for the operations in the above list.
\end{lemma}

\begin{proof}
We just do the case of $\digamma p$.
Suppose $\phi$ is bisimilar to $\phi'$. Let $R$ be the witnessing bisimulation. We define $R'$ between $\digamma p.\phi$
and $\digamma p.\phi'$ simply as $R$ restricted to the vertices not labeled by $p$ in $\phi$. We claim that
$R'$ is a root-preserving bisimulation. Consider any $a,a'$  with $aR'a'$. Consider $ Sa$ in $\digamma p.\phi$.
The only difference with $Sa$ in $\phi$ is that all vertices labeled $p$ are now replaced by the
root with the label of the root. Similarly, for $S'a'$ in  $\digamma p.\phi'$. However, since $R$ was
a root-preserving bisimulation we see that now $R'$ is.
\end{proof}

We will sometimes write $\sigma$, $\tau$, \dots, for substitutions. If the substitution $\sigma$ is
$[q_0:\psi_0,\dots, q_{k-1}:\psi_{k-1}]$ and $\tau$ is $[r_0:\chi_0,\dots, r_{m-1}:\chi_{m-1}]$,
then \[\sigma\star \tau := [q_0:\psi_0,\dots, q_{k-1}:\psi_{k-1},r_0:\chi_0,\dots, r_{m-1}:\chi_{m-1}].\]
This only makes sense if the $q_i$ and the $r_j$ are pairwise disjoint.
 We note an important insight.

\begin{lemma}
\begin{enumerate}[i.]
\item
$\star$ is associative, assuming that the three domains of the substitutions are pairwise disjoint.
\item
Suppose the $q_i$ are disjoint from the domain of $\tau$. Then,
\[\phi[q_0:\psi_0,\dots, q_{k-1}:\psi_{k-1}]\tau \iso \phi([q_0:\psi_0\tau,\dots q_{k-1}:\psi_{k-1}\tau] \star \tau).\]
\end{enumerate} 
\end{lemma}

Suppose $p$ does not occur in $\phi$. We define the following operation:
\begin{itemize}
\item
$\snip{\phi}{p} =:\phi'$ is obtained as follows. 
If the root of $\phi$ is not on a cycle $\phi' :=\phi$. Otherwise, we take an $r^\star$ that is not in $V$, the set of vertices of $\phi$. 
Let 
$V' := V \cup\verz{r^\star}$ and $\lambda' := \lambda \cup \verz{\tupel{r^\star,p}}$. 
We define $f:V \to V'$ by $f(r) = r^\star$ and $f(a)= a$, if $a \neq r$. As usual, we write $f\tupel{a_0,\dots, a_{k-1}}$ for
$\tupel{fa_0,\dots,fa_{n-1}}$.
We take $S'a:= fSa$ if $a\neq r^\star$ and $Sr^\star := \varepsilon$. 
\end{itemize}

\noindent
So, $\snip{\phi}{p}$ is the result of redirecting all incoming arrows of the root, if there are any, to a new vertex labeled $p$.
The new vertex, of course, does not have outgoing arrows. We note that, whether $p$ occurs in $\phi$ or not, $\snip{\phi}{p}$
is modalised in $p$. 

A good heuristic, in case the root is on a cycle, is to view $\snip{\phi}{p}$ as a non-deterministic sub-formula of $\phi$.
The number of nodes increases, so in the sense this `subformula' is larger than
the original formula. On the other hand, trivially, the number of cycles decreases, so in that sense
the `subformula' is smaller. This last feature will be quite useful in the paper.

It will be convenient to write $\snip\phi\psi$ for $\snip \phi p[p:\psi]$.

We enumerate some useful facts about the operation $\snip{\cdot}{\cdot}$.
\begin{lemma}\label{flashssmurf}
Suppose the root of $\phi$ is a cycle vertex and $p$ does not occur in $\phi$. 
Suppose $p$ is modalised in $\psi$. Then,
\begin{enumerate}[i.]
\item
$\mathsf c\snip\phi p< \mathsf c(\phi)$.
\item
$\digamma p.\snip \phi p\iso \phi$.
\item
$\snip \phi \phi \bisim\phi$.
\item
$\digamma p.\psi \bisim \psi[p:\digamma p.\psi]$.
\end{enumerate}
\end{lemma}

\begin{remark}
The operation $\snip{\cdot}{\cdot}$ has to be treated with great care since we may have $\phi\bisim \phi'$ but
$\snip \phi p \not\bisim \snip{\phi'} p$. This can be easily seen from the fact that we can always unravel a
formula a bit to a bisimilar one of which the root is not on a cycle.
\end{remark}

Interestingly, there is something like uniqueness of fixed points modulo bisimulation, which gives
an intriguing analogy with the de Jongh-Sambin-Bernardi Theorem.

\begin{theorem}\label{unibis}
Suppose $p$ is modalised in $\phi$ and $\psi \bisim \phi[p:\psi]$. Then
$\psi \bisim \digamma p.\phi$.
\end{theorem}

\begin{proof}
Let $E$ be the embedding of the vertices of $\phi$ that are not occurrences of $p$ into $\phi[p:\psi]$.
Let $F$ be the embedding of the vertices of $\psi$ into the vertices of the substituted copy of $\psi$ in $\phi[p:\psi]$.
We note that $F$ is a (non-root-preserving) bisimulation between $\psi$ and $\phi[p:\psi]$. 

Let $R$ be the maximal bisimulation between $\phi[p:\psi]$ and  $\psi$.  
We claim that $R^\star := E; R$ is a root-preserving bisimulation
between $\digamma p.\phi$ and $\psi$. (Here `;' is composition in the order of reading.) 

Clearly $R^\star$ connects the roots of $\digamma p.\phi$ and $\psi$ (since the root of $\phi$ cannot be a
$p$-occurrence). 
Consider any node $a$ of  $\digamma p.\phi$. Then $a$ is, by definition, a non-$p$-occurrence
in $\phi$. Suppose $a R^\star b$, say $a E c R b$. Let $a' :=S_{\digamma p.\phi,i}a$ and $b' := S_{\phi[p:\psi],i} b$,
$c' := S_{\psi,i}c$, $d' := S_{\phi,i}a$
\begin{enumerate}[C{a}se 1:]
\item Suppose $d'$ is not a $p$-occurrence in $\phi$.
In this case $a' = d' E c'$. Finally,
 because $cRb$ it follows that $c'Rb'$ and, hence,
$a' R^\star b'$. 
\item Suppose $d'$ is a $p$-occurrence in $\phi$. In this case $a' = r_\phi E r_{\phi[p:\psi]}$. Moreover,
$c'$ is the root of the substituted copy of  $\psi$ in $\phi[p:\psi]$ and $c'R b'$.
We have $r_{\phi[p:\psi]} R r_\psi F c' R b'$. Since, $R;F;R$ is a bisumulation, it is contained in $R$,
so  $r_{\phi[p:\psi]} R b'$. It follows that $a' E r_{\phi[p:\psi]} R b'$, i.o.w., $a'R^\star b'$.
\end{enumerate}
We may conclude that $\digamma p.\phi \bisim \psi$.
\end{proof}

We note that the above proof does not use the full guard condition. It just uses that the root is not a $p$-occurrence.

The fixed point theorem as given in Lemma~\ref{flashssmurf} has a easy generalisation to systems of equations.
Suppose $\mathcal E$ is a system of equations, i.e., a function from a finite set $Q$ of variables to formulas $\phi_q$. 
We form a directed (unlabeled) graph (without ordered successors) ${\sf G}_{\mathcal E}$ with domain $Q$
where we have an arrow from $q$ to $q'$ precisely if $\phi_q$ is not modalised in $q'$.
 We say that $\mathcal E$ is \emph{modalised} if ${\sf G}_{\mathcal E}$ is acyclic. 
 This condition generalises the usual one: $\phi$ is modalised in $p$ iff
the equation $p \mapsto \phi$ is modalised. 

We want to solve $\mathcal E$.
This means that we want to find a function $\mathcal F: q \mapsto \psi_q$ on $Q$,
 such that $\psi_q \bisim \phi_q\mathcal F$, for all $q\in Q$, where in the right-hand-side we view $\mathcal F$ as a substitution.
We demand that the $q'\in Q$ do not occur in the $\psi_q$.
 
We define a new unrooted graph $\Psi$ as follows.
We define a sub-graph ${\sf G}^\circ_{\mathcal E}$ of ${\sf G}_{\mathcal E}$ as follows.
We have an arrow from $q$ to $q'$ iff  $\lambda_{\phi_q}(r_{\phi_q}) = q'$.
This new graph is clearly non-cyclic and every $q$ in $Q$ has at most one outgoing arrow.
Let ${\sf end}(q)$ be the variable at the end of the unique outgoing path from $q$.

We have the following definitions.
\begin{itemize}
\item
$V_\Psi$ is the set of  all pairs $\tupel{\phi_q,a}$, where $q\in Q$ and $a \in V_{\phi_q}$ and $\lambda_{\phi_q}(a) \not\in Q$.
\item
$\lambda_\Psi(\tupel{\phi_q,a}) := \lambda_{\phi_q}(a)$.
\item
Suppose $\lambda_{\phi_q}(a) = q'\in Q$. Then, ${\sf idfy}(\tupel{\phi_q,a}) := \tupel{\phi_{{\sf end}(q')},r_{\phi_{{\sf end}(q')}}}$.
In all other cases, ${\sf idfy}(\tupel{\phi_q,a}) := \tupel{\phi_q,a}$.
We note that the label of a value of {\sf idfy} cannot be in $Q$.
\item
Suppose $S_{\phi_q}a = \tupel{b_0,\dots,b_{n-1}}$. Then,
\[S_\Psi\tupel{\phi_q,a} := \tupel{{\sf idfy}(\tupel{\phi_q,b_0}),\dots, {\sf idfy}(\tupel{\phi_q,b_{n-1}})}.\]
\item
$\mathfrak r_q := {\sf idfy}(\tupel{\phi_q,r_{\phi_q}})$.
In a sense the $\mathfrak r_q$ are multiple roots of $\Psi$.
\item
$\digamma \mathcal E :=\mathcal F$, where $q\mathcal F := \psi_q := \sugr \Psi {\mathfrak r_q}$.
We write $\digamma_{\!\! q}\, \mathcal E$ for $q \digamma \mathcal E$.
\end{itemize}

\begin{theorem}\label{mufipohe}
$\digamma \mathcal E$ is a solution of $\mathcal E$.
\end{theorem}

\begin{proof}
Let $\phi_q$, $\psi_q$ and $\mathcal F$ be as above.
Without loss of generality we may assume that the vertices of
$\Psi$ are disjoint from the vertices of the $\phi_q$. 
We define a bisimulation between $\psi_q$ and $\phi_q\mathcal F$ as follows:
$\alpha R \beta$ iff $\alpha = \beta$ or ($\alpha = \tupel{\phi_q,a}$ and $\beta = a$).

We note that, in case $\phi_q = q'\in Q$, the second disjunct cannot become active.
In this case, $\psi_q$ and $\phi_q\mathcal F$ will be identical.

Suppose that $\phi_q$ is not in $Q$. We consider the case where $\alpha = \tupel{\phi_q,a}$ and $\beta = a$.
We need an auxiliary definition.
 Suppose $\lambda_{\phi_q}(a) = q'\in Q$. Then, ${\sf idfy}_0(a) := \tupel{\phi_{{\sf end}(q')},r_{\phi_{{\sf end}(q')}}}$.
In all other cases, ${\sf idfy}_0(a) := a$. We lift the relation $R$ to sequences in the obvious way.

Let $S_{\phi_q} a = \tupel{b_0,\dots,b_{n-1}}$. 
Then,\qedright
\begin{eqnarray*}
S_{\psi_q}(\alpha) & = &  \tupel{{\sf idfy}(\tupel{\phi_q,b_0}),\dots, {\sf idfy}(\tupel{\phi_q,b_{n-1}})} \\
& R & \tupel{{\sf idfy}_0(b_0),\dots, {\sf idfy}_0(b_{n-1})} \\
& = & S_{\phi_q\mathcal F} a
\end{eqnarray*}
\end{proof}

\section{Cyclic Henkin Logic}\label{cyhelo}

In this section we develop the logic {\sf CHL}. We choose to develop it for itself and not
as part of a wider class of logics, even if, from a systematic standpoint, that would be better.
The reason is simply the desire  not to overburden the presentation. We 
will touch on the broader perspective in Section~\ref{leulosection} and, specifically, in Subsection~\ref{latraint}.

Subsection~\ref{haraldsmurf} provides the basic development of {\sf CHL}. 
Curiously, the central result of the subsection is inter-substitution of equivalents, Theorem~\ref{thec}.
All further development rests on this central result.
In Subsection~\ref{altax}, we consider
some alternative axiomatisations. Finally, in Subsection~\ref{kripsem}, we provide the Kripke semantics for {\sf CHL}.

\subsection{Basic Development}\label{haraldsmurf}
Cyclic Henkin Logic or {\sf CHL} is axiomatised as follows.
\begin{enumerate}[{\sf chl}1.]
\item
If $\vdash \phi$ and $\vdash \phi \to \psi$, then $\vdash \psi$.
\item
If $\vdash \phi$, then $\vdash \opr\phi$.
\item
We have all substitution instances of propositional tautologies.
Here we think of these tautologies as given by the usual parse trees.
\item
$\vdash \opr(\phi \to \psi) \to (\opr \phi \to \opr\psi)$.
\item\label{minismurf}
Suppose $\phi\bisim\psi$. Then $\vdash \phi \iff \psi$.
\item
If $\vdash \opr \phi \to \phi$, then $\vdash \phi$.
\end{enumerate}

We note that if we think of our formulas as the result of dividing out bisimularity, then 
Axiom Scheme {\sf chl}\ref{minismurf} becomes superfluous. That is good since
it is the only principle that does not have the standard form of an axiom scheme.

Our first five axioms amount to something like cyclic {\sf K}. L\"ob's Rule provides this system
with bite. We remind the reader that ordinary {\sf K} is closed under L\"ob's Rule, so it is the presence of circularity that
makes the rule powerful.

We will use the notation $\Gamma \vdash_{\sf CHL} \phi$, where $\Gamma$ is a finite set of formulas
for: ${\sf CHL}\vdash \bigwedge \Gamma \to \phi$. Trivially, we have the deduction rule for this notion.

The following theorem tells us that {\sf CHL} is indeed a logic.

\begin{theorem}\label{thea}
If ${\sf CHL} \vdash \phi$, then ${\sf CHL} \vdash \phi [p_0:\psi_0,\dots,p_{k-1}:\psi_{k-1}]$.
\end{theorem}

\begin{proof}
By a simple induction on proofs. In the case of {\sf chl}5, this uses the safety of substitution.
\end{proof}

\begin{theorem}\label{theb}
\begin{enumerate}[i.]
\item
${\sf CHL}\vdash \phi \iff \snip \phi \phi$.
\item
 ${\sf CHL}\vdash \digamma p.\psi \iff \psi[p:\digamma p.\psi]$, assuming that $\psi$ is modalised in $p$.
 \end{enumerate}
\end{theorem}

\begin{proof}
This is immediate by the fact that $\phi\bisim\snip\phi\phi$ and
$ \digamma p.\psi \bisim \phi[p:\digamma p.\psi]$.
\end{proof}

We define $\opr^\grullet \phi := \digamma p.\opr(\phi \wedge p)$, where $p$ does not occur in $\phi$.

\begin{theorem}\label{stevigesmurf}
We have:
\begin{enumerate}[i.]
\item
If ${\sf CHL} \vdash \opr \phi$, then ${\sf CHL} \vdash \opr^\grullet\phi$.
\item
If ${\sf CHL} \vdash  \phi$, then ${\sf CHL} \vdash \opr^\grullet\phi$. \sbra{Necessitation for $\opr^\grullet$}
\item
If ${\sf CHL} \vdash \opr^\grullet \phi \to \phi$, then ${\sf CHL} \vdash \phi$. \sbra{{\sf LR} for $\opr^\grullet$}
\end{enumerate}
\end{theorem}

\begin{proof}
We verify (i). Suppose ${\sf CHL}\vdash\opr \phi$.
We have:
\begin{eqnarray*}
 \opr\opr^\grullet \phi & \vdash_{\sf CHL} & \opr (\phi \wedge \opr^\grullet \phi) \\
& \vdash_{\sf CHL} & \opr^\grullet \phi
\end{eqnarray*} 
Hence, by L\"ob's rule,  $\vdash_{\sf CHL} \opr^\grullet\phi$. 

(ii) is immediate from (i).

We verify (iii). Suppose $\opr^\grullet\phi \vdash_{\sf CHL} \phi$.
Then, $\opr (\phi \wedge \opr^\grullet\phi) \vdash_{\sf CHL} \phi \wedge \opr^\grullet \phi$.
By {\sf LR}, we find ${\sf CHL}\vdash\phi$.
\end{proof}

\begin{theorem}\label{loebsmurf}
$\opr^\grullet$ satisfies L\"ob's Logic {\sf GL} over ${\sf CHL}$.
\end{theorem} 

\begin{proof}
We have {\sf L1}, i.e. necessitation, by Theorem~\ref{stevigesmurf}(ii).
We verify {\sf L2}.{\small
\begin{eqnarray*}
 \opr(\opr^\grullet (\phi \to \psi) \to (\opr^\grullet \phi \to \opr^\grullet \psi)) \vdash_{\sf CHL}\;\; 
 \opr^\grullet (\phi \to \psi) & \to & \opr((\phi \to \psi) \wedge \opr^\grullet(\phi \to \psi))\\
 & \to & \opr((\phi \to \psi) \wedge (\opr^\grullet \phi \to \opr^\grullet\psi))\\
 & \to & \opr ((\phi \wedge\opr^\grullet\phi) \to (\psi\wedge\opr^\grullet\psi))\\ 
 & \to & (\opr(\phi \wedge \opr^\grullet \phi) \to \opr(\psi \wedge \opr^\grullet \psi)) \\
 & \to & (\opr^\grullet \phi \to \opr^\grullet\psi)  
\end{eqnarray*}
}
So, by L\"ob's Rule, we are done.
We verify {\sf L3}.
\begin{eqnarray*}
\opr(\opr^\grullet\phi \to \opr^\grullet\opr^\grullet\phi) \vdash_{\sf CHL}\;\; \opr^\grullet \phi & \to & \opr\opr^\grullet \phi  \\
& \to & \opr(\opr^\grullet\phi \wedge \opr^\grullet\opr^\grullet \phi) \\
& \to & \opr^\grullet\opr^\grullet\phi
\end{eqnarray*}
By L\"ob's Rule we are done. We verify L\"ob's Principle {\sf L4}. We have:
\begin{eqnarray*}
\opr^\grullet (\opr^\grullet (\opr^\grullet \phi \to \phi) \to \opr^\grullet \phi) \vdash_{\sf CHL}\;\;  
\opr^\grullet (\opr^\grullet \phi \to \phi) & \to & \opr^\grullet\opr^\grullet (\opr^\grullet \phi \to \phi)\\
& \to & \opr^\grullet\opr^\grullet \phi \\
& \to & \opr^\grullet\phi
\end{eqnarray*}
By {\sf LR} for $\opr^\grullet$ (Theorem~\ref{stevigesmurf}(iii)), we are done.
\end{proof}

We write $\dotbox \phi$ for $\phi \wedge \opr\phi$ and, similarly for $\opr^\grullet$.

\begin{corollary}\label{cora}
We have strengthened L\"ob's Rule for $\opr^\grullet$ over ${\sf CHL}$, i.e.,
\[\text{ if }\bigwedge_{i<n}\opr^\grullet \psi_i ,\bigwedge_{j<m}\dotbox^\grullet \chi_j \vdash \opr\phi\to\phi,\text{ then }
\bigwedge_{i<n}\opr^\grullet \psi_i ,\bigwedge_{j<m}\dotbox^\grullet \chi_j \vdash \phi.\]
\end{corollary}

We may simply copy the usual derivation
of strengthened L\"ob's Rule from L\"ob's Principle. For completeness we repeat the argument.

\begin{proof}
Suppose  $\bigwedge_{i<n}\opr^\grullet \psi_i ,\bigwedge_{j<m}\dotbox^\grullet \chi_j \vdash \opr\phi\to\phi$.
Then,  \[(\dag) \;\;\; \bigwedge_{i<n}\opr^\grullet \psi_i ,\bigwedge_{j<m}\dotbox^\grullet \chi_j \vdash \opr^\grullet\phi\to\phi.\]
It follows that $\bigwedge_{i<n}\opr^\grullet \psi_i ,\bigwedge_{j<m}\dotbox^\grullet \chi_j \vdash \opr^\grullet( \opr^\grullet\phi\to\phi)$ and, hence, that
\[(\ddag) \;\;\; \bigwedge_{i<n}\opr^\grullet \psi_i ,\bigwedge_{j<m}\dotbox^\grullet \chi_j \vdash \opr^\grullet \phi.\]
Combining (\dag) and (\ddag), we find: $\bigwedge_{i<n}\opr^\grullet \psi_i ,\bigwedge_{j<m}\dotbox^\grullet \chi_j \vdash \phi$
\end{proof}

Now we are ready to prove a principle that will be a central tool: a strong form of
substitution of equivalents.

\begin{theorem}\label{thec}
Suppose that the variables $s_i$, for $i< n$ and $r_j$ for $j<m$ are pairwise distinct;
  and that $\phi$ is modalised in the $r_i$.
  Then, we have:
\[ \bigwedge_{i<n}  \dotbox^\grullet(\psi_i \iff \chi_i), \bigwedge_{j<m}\opr^\grullet (\theta_j \iff \rho_j) \vdash_{\sf CHL}
 \phi\subun{\vec s:= \vec \psi, \vec r:=\vec \theta\,}
\iff  \phi\subun{\vec s:= \vec \chi, \vec r:=\vec \rho\,}.\]
\end{theorem}

\begin{proof}
Suppose a formula $\psi$ and substitutions $\sigma$ and $\tau$ on $q_0,\dots, q_{k-1}$ are given. 
We consider a conjunction $\alpha =\bigwedge_{i<k} \apr_i(q_i\sigma \iff q_i\tau)$, where $\apr_i$
is either $\dotbox^\grullet$ or $\opr^\grullet$. We say that $\alpha$ is \emph{acceptable}
for $\psi,\sigma,\tau$ if, whenever $\apr_i$ is $\opr^\grullet$, then $\psi$ is modalised in $q_i$. 

We prove: for every $\phi$, we have $P(\phi)$, where:
\begin{eqnarray*}
P(\phi): &\Iff &\text{for every $\vec q$, for every $\sigma,\tau$ on $\vec q$,}\\
&& \text{for every $\alpha$ acceptable 
for $\phi,\sigma,\tau$, we have $\alpha \vdash \phi\sigma \iff \phi\tau$.}
\end{eqnarray*}

The proof is by course-of-values induction  for $\mathsf c(\phi)$.
As a sub-induction we use guard-induction using ${\sf bo}^\circ(\phi)$, the set
of cycle $\opr$-occurrences, as a guard. So, in the sub-induction we prove the induction
property for $\sugr{\phi}{a}$.

We treat the examples of conjunction and box, spitting the second case into the sub-cases that
the occurrence is not on a cycle or is on a cycle.

Suppose $a$ is labeled $\wedge$. Suppose $\alpha$ is acceptable for $(\sugr\phi a),\sigma,\tau$.
Clearly, $\alpha$ is also acceptable for the $(\sugr\phi {S_ia}),\sigma,\tau$.
We have, by the safety of substitution, $\sugr \phi a\bisim(\sugr \phi {S_0a}\wedge \sugr \phi{S_1a})$ and, hence,
$(\sugr \phi a)\sigma\bisim((\sugr \phi {S_0a})\sigma\wedge (\sugr \phi {S_1a})\sigma)$ and, similarly, for $\tau$.
So, assuming the desired propery for the $\sugr \phi {S_ia}$, we are done.

Suppose $a$ is labeled with a box, but not on a cycle.
 We have $\sugr \phi a\bisim  \opr(\sugr\phi S_0a)$. 
 
Let $\alpha$ be acceptable for $(\sugr\phi a),\sigma,\tau$.
Let $\alpha^\star :=\bigwedge_{i<k} \dotbox^\grullet(q_i\sigma \iff q_i\tau)$.
Clearly, $\alpha^\star$ is acceptable for $(\sugr \phi S_0a),\sigma,\tau$.
 By guard-induction, we may assume:
   $\alpha^\star \vdash_{\sf CHL} (\sugr\phi b)\sigma \iff (\sugr\phi b)\tau$.
It follows, by {\sf K}-reasoning, that 
\[ \opr \alpha^\star \vdash_{\sf CHL} \opr((\sugr \phi b)\sigma)\iff  \opr((\sugr\phi b)\tau).\]
 Moreover, $\alpha \vdash_{\sf CHL} \opr\alpha^\star$.
So, 
$\alpha \vdash_{\sf CHL} (\opr(\sugr \phi b))\sigma\iff  (\opr(\sugr \phi b))\tau$, as desired.

Suppose  $a$ is a box-occurrence on a cycle, i.e., $a$ is in the chosen guard. 
Let $\phi_0$ be $\sugr\phi a$. Suppose $\alpha$ is acceptable for $\phi_0,\sigma,\tau$.
We choose $p$ not in $\phi_0$, nor in the $\psi_i$, $\chi_i$, $\theta_j$ and $\rho_j$
and distinct from the $q_i$ and $r_j$.
We have:
\begin{eqnarray*}
\phi_0 \sigma & \bisim & \snip{\phi_0} p[p:\phi_0] \sigma \\
& \bisim & \snip{\phi_0} p (\sigma\star[p:\phi_0\sigma]). 
\end{eqnarray*}
Similarly,  $\phi_0\tau\bisim  \snip{\phi_0} p (\tau \star [p:\phi_0\tau])$.

Since $\mathsf c\snip {\phi_0} p < \mathsf c(\phi_0) \leq \mathsf c(\phi)$, we may apply the main induction hypothesis.
Since $\snip{\phi_0}p$ is modalised in $p$, we find that $\alpha\wedge  \opr^\grullet (\phi_0 \sigma \iff \phi_0\tau)$ is acceptable
for $\snip{\phi_0}p, \sigma \star[p:\phi_0\sigma],  \tau \star[p:\phi_0\tau]$.
Thus, by the main induction hypothesis:
\begin{eqnarray*}
\alpha, \opr^\grullet (\phi_0 \sigma \iff \phi_0\tau)  \vdash_{\sf CHL} 
\phi_0 \sigma & \iff & \snip{\phi_0} p (\sigma\star[p:\phi_0\sigma]) \\
& \iff & \snip{\phi_0} p (\tau\star[p:\phi_0\tau])\\
& \iff & \phi_0\tau
\end{eqnarray*}
By the strengthend L\"ob's Rule, we may omit the assumption $\opr^\grullet (\phi_0 \sigma \iff \phi_0\tau)$ and we are done.
\end{proof}

\begin{theorem}\label{equismurf}
Suppose $\phi$ and $\psi$ are modalised in $p$ and
 ${\sf CHL} \vdash \phi \iff  \psi$, then ${\sf CHL} \vdash  \digamma p.\phi \iff \digamma p.\psi$. 
\end{theorem}

\begin{proof}
Suppose ${\sf CHL} \vdash \phi  \iff  \psi $. Then, 
\begin{eqnarray*}
\opr^\grullet(\digamma p.\phi \iff \digamma p.\psi) \vdash_{\sf CHL}   \digamma p.\phi &  \iff & \phi [p:\digamma p. \phi] \\
& \iff &  \psi  [p:\digamma p.\phi] \\
& \iff &  \psi [p:\digamma p.\psi] \\
& \iff & \digamma p.\psi
\end{eqnarray*}
So, by L\"ob's rule, we are done.
\end{proof}

We  prove the de Jong-Sambin-Bernardi Theorem about the uniqueness of fixed points.

\begin{theorem}
Suppose $\phi$ is modalised in $p$. Then \[{\sf CHL} \vdash \dotbox^\grullet(p \iff \phi) \to (p \iff \digamma p.\phi).\]
\end{theorem}

\begin{proof}
We have:
\begin{eqnarray*}
\dotbox^\grullet(p \iff \phi)\vdash_{\sf CHL} \opr^\grullet (p\iff \digamma p.\phi)  \to  (\;p &\iff & \phi \\
 &\iff & \phi[p:\digamma p.\phi] \\
&\iff & \digamma p.\phi \;)
\end{eqnarray*}
By the Strengthened L\"ob's Rule, we may omit the assumption $ \opr^\grullet (p\iff \digamma p.\phi) $.
\end{proof}

We generalise the de Jong-Sambin-Bernardi Theorem to systems of equations as follows.

\begin{theorem}\label{gejosam}
Suppose $\mathcal E$ is modalised. Then,
\[\bigwedge_{q\in Q}\dotbox^\grullet (q\iff q\mathcal E) \vdash_{\sf CHL} \bigwedge_{q\in Q} (q\iff \digamma_{\!\! q}\,\mathcal E).\] 
\end{theorem}

\begin{proof}
We write 
$\phi_q =: q\mathcal E$ and $\psi_q =:  \digamma_{\!\! q}\,\mathcal E$ and $\mathcal F := \digamma \mathcal E$.
Let
 $\chi := \bigwedge_{q\in Q}\dotbox^\grullet (q\iff \phi_q)$ and $\rho := \opr^\grullet \bigwedge_{q\in Q}(q\iff \psi_q)$. 
We prove $\chi, \rho \vdash q \iff \psi_q$, by induction on ${\sf G}_\mathcal E$. 

Let $Q_q$ be the set of all $q'$ that can be reached from $q$ in ${\sf G}_{\mathcal E}$ via a non-empty path.
We suppose we have the desired result for all $q'$ in $Q_q$. 
Let $\mathcal F_q$ be restriction of $\mathcal F$ to the $Q_q$.
It follows that we have:  
\begin{eqnarray*}
\chi,\rho \; \vdash_{\sf CHL}\;  q &  \iff  & \phi_q \\
& \iff & \phi_q \mathcal F_q \\
& \iff &  \phi_q \mathcal F_q\mathcal F \\
& \iff & \phi_q \mathcal F \\
& \iff & \psi_q 
\end{eqnarray*}
The first equivalence is by $\chi$. The second equivalence follows by the combination of $\rho$ and the
induction hypothesis. The third equivalence follows by $\rho$ in combination with the fact that
all variables from $Q$ in $ \phi_q \mathcal F_q$ are guarded: the only variables that $\phi_q$ can `see'
in ${\sf G}_{\mathcal E}$ have been removed by the substitution. Moreover, no variables from $Q$ occur
in the $\psi_{q'}$ that are substituted. The fourth and fifth equivalence are immediate.

We have shown that:
\[ \chi, \opr^\grullet \bigwedge_{q\in Q}(q\iff \psi_q) \vdash_{\sf CHL} \bigwedge_{q\in Q}(q\iff \psi_q). \]
So, by the Strengthed L\"ob's Rule, we are done.
\end{proof}

Consider a formula $\phi$. We assign to each $\opr$-occurrence $a$ a propositional variable $q_a$, where the
$q_a$ are pairwise distinct and also distinct from the propositional variables of $\phi$. We map the
nodes of $\phi$ to formulas of the language of ordinary modal logic, i.e., the cycle-free formulas, as follows:
\begin{itemize}
\item
${\sf E}_\phi(a) := q_a$, if $a$ is a $\opr$-occurrence.
\item
${\sf E}_\phi(a) := ({\sf E}_\phi(S_0a) \wedge {\sf E}_\phi(S_1a))$, if $a$ is a $\wedge$-occurrence.
Similarly, for the other connectives and for the propositional variable-occurrences in $\phi$.
\end{itemize}

Our definition is correct by guard-recursion. 
We write $\psi_a$ for ${\sf E}_\phi(a)$. 

\begin{theorem}\label{thed}
Let $G$ be the set of $\opr$-occurrences of $\phi$. We have:
\[\bigwedge_{a\in G} \dotbox^\grullet(q_a\iff \opr \psi_{S_0a}) \vdash_{\sf CHL} \psi_r \iff \phi.\]
\end{theorem}

\begin{proof}
We show by guard induction with ${\sf bo}(\phi)$ as guard that, for all nodes $b$, we have:
\[ \opr^\grullet\bigwedge_{e\in V}(\psi_e\iff \sugr\phi e),
\bigwedge_{a\in G} \dotbox^\grullet(q_a\iff \opr \psi_{S_0a}) \vdash \psi_b \iff \sugr \phi b.\]
We first treat the case of conjunction. Suppose $b$ is a $\wedge$-occurrence.
We assume our desired conclusion for $S_0b$ and $S_1b$.
We have:
\begin{eqnarray*}
 \opr^\grullet\bigwedge_{e\in V}(\psi_e\iff \sugr\phi e), 
\bigwedge_{a\in G} \dotbox^\grullet(q_a\iff \opr \psi_{S_0a}) \vdash_{\sf CHL}\;\;  \psi_b & \iff & (\psi_{S_0b} \wedge \psi_{S_1b}) \\
& \iff & ((\sugr \phi {S_0b}) \wedge (\sugr \phi {S_1b})) \\
& \iff & \sugr \phi b
\end{eqnarray*}
The reasoning for the propositional variable-occurrences of $\phi$ and for the other non-box connectives is similar.
Suppose $b$ is a $\opr$-occurrence. We have:
\begin{eqnarray*}
\opr^\grullet\bigwedge_{e\in V} (\psi_e\iff \sugr\phi e),
\bigwedge_{a\in G} \dotbox^\grullet(q_a\iff \opr \psi_{S_0a}) \vdash_{\sf CHL}\;\; q_b & \iff & \opr\psi_{S_0b} \\
& \iff  & \opr (\sugr\phi {S_0b}) \\
& \iff  & \sugr \phi b
\end{eqnarray*}
So we find:
\[ \opr^\grullet\bigwedge_{e\in V} (\psi(e)\iff \sugr\phi e),
\bigwedge_{a\in G} \dotbox^\grullet(q_a\iff \opr \psi_{S_0a}) \vdash_{\sf CHL} \bigwedge_{e\in V} (\psi_e\iff \sugr\phi e) \]
We apply the strengthened version of L\"ob's Rule to obtain:
\[ 
\bigwedge_{a\in G} \dotbox^\grullet(q_a\iff \opr \psi_{S_0a}) \vdash_{\sf CHL} \bigwedge_{e\in V} (\psi_e\iff \sugr\phi e) \]
For this last insight, the desired conclusion is immediate.
\end{proof}

\begin{remark}
One would hope that the same kind of treatment we give here for {\sf CHL} would be also possible for the
$\mu$-calculus. Formulas are defined in the same way, only the guarding constraint is replaced by the following constraint.
Consider any directive cycle $\gamma:C \to \omega$ in $\phi$. Consider the set $X$ consisting  of all occurrences $a$ of $\neg$ in $C$
plus all occurrences $b$ of $\to$ in $C$ such that $\gamma(b)=0$. We demand that $X$ has an even number of elements.
In stead of L\"ob's Rule one would have the Minimality Rule: if $\vdash \snip \phi \alpha \to \alpha$, then $\vdash \phi \to \alpha$.
I have no idea how the details of this will work out. What replaces guard-induction and recursion?
\end{remark}

\subsection{Alternative Axiomatisations}\label{altax}
We provide some alternative axiomatisations for {\sf CHL}.
We consider the following axioms an rules.
\begin{description}
\item[{\sf IPE}]
 Suppose $\phi$ and $\psi$ are modalised in $p$ and $\vdash \phi \iff \psi$, then $\vdash \digamma p.\phi \iff \digamma p.\psi$
 (Intersubstitutivity of Provable Equivalents).
 \item[${\sf N}^\grullet$]
 If $\vdash \phi$, then $\vdash \opr^\grullet \phi$.
 \item[{\sf H}]
 $\vdash {\sf H}$, where ${\sf H}= \digamma p.\opr p$.
 \end{description}
 
 We define the following theories:
 \begin{itemize}
 \item
 ${\sf K}^{\circ-}$ is given by {\sf chl}1-5.
 \item
 ${\sf K}^\circ$ is $K^{\circ-}$ plus {\sf IPE}.
 \item
 ${\sf CHL}_0$ is $K^{\circ}$ plus {\sf H}.
 \item
 ${\sf CHL}_1$ is ${\sf K}^{\circ -}$ plus ${\sf N}^\grullet$.
 \end{itemize}
 
 Here ${\sf K}^\circ$ is the reasonable circular version of {\sf K}.
 We will show that over this theory L\"ob's rule is equivalent with the
 Henkin sentence (as axiom).
 
 \begin{theorem}
  ${\sf CHL}_0$ and  ${\sf CHL}_1$ prove the same theorems as {\sf CHL}.
 \end{theorem}
 
 \begin{proof}
 To prove inclusion of ${\sf CHL}_1$ in ${\sf CHL}_0$, we show that ${\sf CHL}_0$ is closed under
 ${\sf N}^\grullet$ Suppose ${\sf CHL}_0 \vdash \phi$. It follows that ${\sf CHL}_0\vdash \opr (p\wedge \phi) \iff \opr p$.
 Hence, by {\sf IPE}, we have  ${\sf CHL}_0\vdash \digamma p. \opr (\phi \wedge p) \iff\digamma p. \opr p$,
 i.e., ${\sf CHL}_0\vdash \opr^\grullet  \phi \iff {\sf H}$. So, by axiom {\sf H}, we find ${\sf CHL}_0 \vdash \opr^\grullet\phi$.
 
 To prove inclusion of {\sf CHL} in ${\sf CHL}_1$ it is sufficient to show that  ${\sf CHL}_1$ is closed under L\"ob's Rule.
 This uses a well-known proof of L\"ob's Rule.
 Suppose ${\sf CHL}_1\vdash \opr\phi \to \phi$. Let $\nu := \digamma q.\opr^\grullet (q\to \phi)$, where $q$ does not
 occur in $\phi$. We have:
 \begin{eqnarray*} 
 {\sf CHL}_1 \vdash \nu & \to & \opr^\grullet (\nu \to \phi) \\
 & \to & \opr( (\nu \to \phi) \wedge \opr^\grullet(\nu \to \phi)) \\
 & \to & \opr((\nu \to \phi) \wedge \nu) \\
 & \to & \opr \phi \\
 & \to & \phi
 \end{eqnarray*}
 So, we have (a)  ${\sf CHL}_1\vdash \nu \to \phi$. Hence, by ${\sf N}^\grullet$, we find ${\sf CHL}_1\vdash \opr^\grullet(\nu \to \phi)$.
 It follows that (b)  ${\sf CHL}_1\vdash \nu$. Combining (a) and (b), we find: $\vdash \phi$.
 
 The inclusion of ${\sf CHL}_0$ in {\sf CHL} follows from our previous results.
 \end{proof}
 
 \subsection{Kripke Semantics}\label{kripsem}
A Kripke model for {\sf CHL} is given by a triple $\tupel{W,\sqsubset,f}$, where $W$ is a non-empty set of worlds,
$\sqsubset$ is an acyclic binary relation on $W$ and $f: W \times {\sf Prop} \to \verz{0,1}$, where {\sf Prop} is the set of
propositional variables.

Consider any formula $\phi$. We define ${\sf Ev}_\phi:W \times V \to \verz{0,1}$ as follows.
We define ${\sf Ev}_\phi(w,\cdot)$ assuming that we already have defined ${\sf Ev}_\phi(w',\cdot)$
for all $w'\sqsupset w$. We use guard-recursion w.r.t. ${\sf bo}(\phi)$ as guard.
\begin{itemize}
\item
If $a$ is an occurrence of $p$, then ${\sf Ev}_\phi(w,a) = f(w,p)$.
\item
If $a$ is an occurrence of $\wedge$, then
\[{\sf Ev}_\phi(w,a) = {\sf min}({\sf Ev}_\phi(w,S_0a),{\sf Ev}_\phi(w,S_1a)),\] and, similarly for the other non-box
connectives.
\item
If $a$ is a $\opr$-occurrence, we set
${\sf Ev}_\phi(w,a) =1$ iff, for all $w'\sqsupset w$, we have ${\sf Ev}_\phi(w',S_0a)=1$.
We take ${\sf Ev}_\phi(w,a)=0$, otherwise.
\end{itemize}
We note that this definition works since in the clause for box do not call on values for $w$ but on the
previously defined values for $w'\sqsupset w$.

We define $\eva\phi (w) := {\sf Ev}_\phi(w,r_\phi)$ and we write $w \Vdash \phi$ for $\eva \phi (w)=1$.
Since our guard is preserved by bisimulation, we find that 
if $\phi \bisim \phi'$, then $\eva \phi(w) = \eva {\phi'} (w)$.

We easily verify the usual clauses like $w \Vdash \phi\wedge \psi$ iff $w\Vdash \phi$ and $w\Vdash \psi$.
Using this the validity of {\sf CHL} is immediate.

\begin{remark}
We can derive the Kripke completeness theorem for {\sf CHL} if finite acyclic models in two ways. The first is using the synonymy of
{\sf CHL} and the Henkin Calculus {\sf HC}, which is essentially {\sf CHL} on standard syntax using a 
fixed point operator. We can prove the Kripke completeness for {\sf HC}, for example, by showing
that {\sf HC} is synonymous to the well-founded part of the $\mu$-calculus and invoking the Kripke compleness theorem
for the $\mu$-calculus. We will give the details of this argument in a sequel to this paper by Tadeusz Litak
and myself. Alternatively, we can use the Kripke Completeness of {\sf WfL} and use the reduction of
{\sf CHL} to that theory: see Subsection~\ref{geaardesmurf} for more detail. 
Of course, it would be much nicer to give a direct proof
of the completeness theorem in terms of the circular syntax. To do this remains open.
\end{remark}
 
 \section{L\"ob's Logic meets Cyclic Syntax}\label{leulosection}
 The counterpart of {\sf CHL} in cycle-free language is the Henkin Calculus {\sf HC}. This calculus employs a variable-binding
 fixed point operator in the object language. The logics {\sf CHL} and {\sf HC} are synonymous/definitionally equivalent.
 We will establish this fact in a later paper in which {\sf HC} is developed. Here we will treat the simpler
 case of the synonymy of ${\sf GL}^\circ$, the extension of {\sf CHL} with the transitivity axiom scheme {\sf L3}, aka 4, and
 L\"ob's Logic {\sf GL}.
 
 In Subsection~\ref{latraint}, we set up a modest framework for a province in which both theories
 live and some basics of comparing theories living in that province. In Subsection~\ref{lolocyc}, the logic
 ${\sf GL}^\circ$ is introduced and in Subsection~\ref{lowicyc} we do the same for ordinary {\sf GL} (as it appears
 in our framework). In Subsection~\ref{josaint}, we develop the de Jongh-Sambin algorithm as one half of the witness
 of a synonymy between ${\sf GL}^\circ$ and {\sf GL}.
 
 \subsection{Languages, Translations and Interpretations}\label{latraint}
 In this section we look at a restricted class of logics and employ a very restricted framework
 of interpretations. The interpretations are something like \emph{$\opr$-preserving logic-interpretations}.
 However, in this section we will call them \emph{$\ell$-interpretations}. 
 
 A \emph{language} $\mathcal L$ in this section will simply be a sub-set of our full language
 $\mathbb L^\circ$ that is closed under (i) bisimilarity, (ii) the propositional variables
 and the syntactic operations associated with the connectives, (iii) subformulas and (iv) substitution.
 The minimal language is $\mathbb L$, the set of all acyclic formulas. 
 We may view $\mathbb L$ as the ordinary language of modal logic since 
 each bisimulation equivalence class of an acyclic formula has a unique finite tree (modulo isomorphism)
 as normal form.
 
 We define ${\sf HL}(\mathcal L)$ as the logic in the language $\mathcal L$ that is axiomatised by 
 {\sf chl}1--6. One easily sees that this definition makes sense. We note that 
 ${\sf CHL} = {\sf HL}(\mathbb L^\circ)$ and that
 {\sf K} proves the same theorems as ${\sf HL}(\mathbb L)$.
 A \emph{logic} $\Lambda$ will be a rule-preserving extension of one of the  ${\sf HL}(\mathcal L)$ in the same language
 by schematic rules.
 
 A salient language is $\mathbb L^\grullet$. This
 is the language generated by the variables, the logical connectives including $\opr$ and $\opr^\grullet$.
We define ${\sf WfL} := {\sf HL}(\mathbb L^\grullet)$.

Here is an important observation.
\begin{observation}
Suppose $\mathcal L$ is a language that extends $\mathbb L^\grullet$. Then
Theorems~\ref{thea}, \ref{theb}\sbra i, \ref{stevigesmurf}, \ref{loebsmurf}, \ref{thec} and
Corollary~\ref{cora} remain valid when we replace {\sf CHL} by ${\sf HL}(\mathcal L)$.
 \end{observation}

We can simply check that these results do not use  closure under the operation $\digamma p.(\cdot)$.
  
 A \emph{local translation} $\mathscr T$ of a formula $\phi$ into a logic $\Lambda$
 is a mapping from $V_\phi$ to the language of $\Lambda$, such that  
 if $a$ is $\$$-occurrence for  a connective $\$$ of arity $n$ or a variable treated as having arity 0, where
 \[\Lambda \vdash {\mathscr T}(a) \iff \$({\mathscr T}(S_0a),\dots, {\mathscr T}(S_{n-1}a)).\]
 Strictly speaking the local translation is given as the triple of the formula, the function and the logic. 
 
\begin{theorem}\label{unitrans}
Consider any logic $\Lambda$. Suppose $\mathscr T$ is a local translation of
of $\phi$ in $\Lambda$ and $\mathscr T'$ is a local translation of
of $\phi'$ in $\Lambda$. Suppose further that 
 $R$ is a bisimulation between $\phi$ and $\phi'$, where $R$ does not need to be root-preserving.
Suppose $aRa'$. Then, $\Lambda \vdash \mathscr T(a) \iff \mathscr T'(a')$.
\end{theorem}

\begin{proof}
Let $\chi := \bigwedge_{bRb'}  (\mathscr T(b) \iff \mathscr T'(b'))$.
We prove by guard-induction on $\opr$-occurrences $d$ in $\phi$ that, if $dRd'$, then
 $\opr\chi \vdash_\Lambda  \mathscr T(d) \iff \mathscr T'(d')$. The cases of the non-box occurrences
 are trivial. In case $d$ is a $\opr$-occurrence, we have:
 \begin{eqnarray*}
  \opr\chi \vdash  \mathscr T(d) & \iff & \opr  \mathscr T(S_0d) \\
  & \iff & \opr \mathscr T'(S'_0d') \\
  & \iff &  \mathscr T'(d')
  \end{eqnarray*}
  It follows that $\opr\chi \vdash_\Lambda \chi$. So, by L\"ob's rule, we are done.
\end{proof}
 
 A \emph{global translation} $\mathscr G$ of a language $\mathcal L$ into a logic $\Lambda$ is
 a function from $\mathcal L$ to the language of $\Lambda$ that commutes modulo
 $\Lambda$-provable equivalence with the propositional variables and the connectives.
 Strictly speaking a global translation is given as the triple of language, function and logic.
 We usually omit the `global' of `global translation'.
 
 We collect some trivial insights.
 
 \begin{theorem}
 \begin{enumerate}[i.]
 \item
 Suppose $\mathscr G$ is a global translation of $\mathcal L$ into  $\Lambda$. Let $\phi$ be
 an $\mathcal L$-formula. For $a\in V_\phi$, we define $\mathscr T(a) := \mathscr G(\sugr\phi a)$.
 Then, $\mathscr T$ is a local translation of $\phi$ in $\Lambda$.
 \item
 Suppose every $\phi$ in $\mathcal L$ has a local translation $\mathscr T_\phi$ in $\Lambda$.
 Then  $\mathscr G$ with $\mathscr G(\phi) := \mathscr T_\phi(r_\phi)$ is a global translation of
 $\mathcal L$ in $\Lambda$.
 \item
 If $\mathscr G$ and $\mathscr G'$ are global translations of $\mathcal L$ in $\Lambda$.
 Then, for each $\phi$, we have $\Lambda \vdash \mathscr G(\phi) \iff \mathscr G'(\phi)$.
 \end{enumerate}
 \end{theorem}
 
 A translation $\mathscr G$ is an $\ell$-translation if it commutes with substitution.
 This means that, for all formulas $\phi$ and all $\mathcal L$-substitutions $\sigma$, 
 we have $\Lambda \vdash \mathscr G(\phi\sigma) \iff \mathscr G(\phi) (\mathscr G \circ\sigma)$. 
 
 A interpretation $\mathscr K:\Lambda \to \Lambda'$ is given as a triple of $\Lambda$, a function $\mathscr G$,
 and $\Lambda'$, where $\mathscr G$ is a translation from the language of $\Lambda$ into $\Lambda'$.
 We demand that if $\Lambda \vdash \phi$, then $\Lambda' \vdash \mathscr G(\phi)$.
 Usually, we will confuse $\mathscr K$ with its underlying translation, writing e.g. $\mathscr K(\phi)$.
 
 An interpretation is an $\ell$-interpretation iff the underlying translation is an $\ell$-translation.  
 
 We count two interpretations as the same if their values are provably equivalent in the target logic.
 The identity interpretation ${\sf ID}_\Lambda$ on $\Lambda$ is the interpretation based on the identity 
 function ${\sf Id}_{\mathcal L}$, where $\mathcal L$ is the language of $\Lambda$.
 We define composition of interpretations in the expected way. We leave the simple 
 verification that interpretations are closed under composition to the reader.
 It is also easy to see that the composition of $\ell$-interpretations delivers a $\ell$-interpretation.
 
 We note that the fact that we have at most one interpretation (from our restricted class) between two
 logics implies that, whenever we have mutual interpretability, we have synonymy: the composition of
 an interpretation from $\Lambda$ to $\Lambda'$ an an interpretation from $\Lambda'$ to $\Lambda$
 will be an interpretation from $\Lambda$ to $\Lambda$ and this must be the identity.
 
 \subsection{L\"ob's Logic with Cycles}\label{lolocyc}
 We define ${\sf GL}^\circ$ as {\sf CHL} plus {\sf L3} or 4, to wit $\vdash \opr\phi \to \opr\opr\phi$.
 Here the intended notion of extension is one that preserves the rules.
 
The theory ${\sf GL}^\circ$ is extensionally the same as the theory ${\sf K}^{\circ -}$ plus {\sf L3}.
This is because in the presence of the fixed points we may prove
L\"ob's Principle in the usual way. Then, the desired closure under L\"ob's rule follows, by the usual proof
of L\"ob's Rule from L\"ob's Principle. We will write $\equiv^\circ$ for provable equivalence in ${\sf GL}^\circ$.

We show that in ${\sf GL}^\circ$ the modalities $\opr$ and $\opr^\grullet$ coincide.

\begin{theorem}
${\sf GL}^\circ \vdash \opr^\grullet \phi \iff \opr \phi$.
\end{theorem}

\begin{proof}
We have: \qedright
\begin{eqnarray*}
\opr (\opr^\grullet \phi \iff \opr \phi)\; \vdash_{{\sf GL}^\circ} \; \opr^\grullet\phi & \iff & \opr(\phi \wedge \opr^\grullet\phi) \\
& \iff & \opr(\phi \wedge \opr\phi) \\
& \iff & \opr\phi
\end{eqnarray*}
\end{proof}

Thus, all insights that we accumulated for $\opr^\grullet$ in {\sf CHL} transfer to $\opr$ in ${\sf GL}^\circ$.

Here is a ${\sf GL}^\circ$ reformulation of a well-known insight due to Dick de Jongh. 

\begin{lemma}\label{kleutersmurf}
Suppose the root of $\psi$ is an $\opr$-occurrence.
Then, we have $\psi \equiv^\circ \snip \psi \top$. 
\end{lemma}

\begin{proof}
If the root of $\psi$ is not on a cycle, this is trivial. So, suppose the root is on a cycle.
Let $\chi := {\sf su}_0(\psi)$. We have $\psi \equiv^\circ \opr \chi$, and hence
$\glc \vdash \psi \to \opr\psi$. Similarly, for $\snip \psi \top$.
Suppose $p$ does not occur in $\psi$. We have:
\begin{eqnarray*}
\psi  & \vdash_{\glc} & \opr \psi \wedge \snip \psi p[p:\psi] \\
 & \vdash_{\glc} & \opr (\psi \iff \top) \wedge \snip \psi p[p:\psi] \\
  & \vdash_{\glc} &  \snip \psi p[p:\top] \\
   & \vdash_{\glc} &  \snip \psi \top 
\end{eqnarray*}

We also have:
\begin{eqnarray*}
\opr\psi, \snip\psi \top  & \vdash_{\glc} & \opr (\psi \iff \top)  \wedge \snip \psi p[p:\top] \\
 & \vdash_{\glc} &  \snip \psi p[p:\psi] \\
  & \vdash_{\glc} & \psi 
  \end{eqnarray*}
  By the Strengthed L\"ob's rule, we find $\snip \psi \top \vdash_{\glc}  \psi$. 
\end{proof}

\subsection{L\"ob's Logic without Cycles}\label{lowicyc}

In our context we may define L\"ob's Logic {\sf GL} simply as
${\sf HL}(\mathbb L)$ plus {\sf L}3. We note that this is
{\sf K} plus L\"ob's Rule plus {\sf L}3.
We will employ the usual facts about {\sf GL} and especially the
de Jongh result:

\begin{lemma}\label{dejosmurf}
${\sf GL} \vdash (\opr\phi)[p:\top] \iff (\opr\phi) [p:(\opr\phi)[p:\top]]$.
\end{lemma}

\subsection{The de Jongh-Sambin Interpretation}\label{josaint}
We define functions {\sf js} and ${\sf js}^\star$. Our aim is to show that ${\sf js}^\star(\phi,\cdot)$ is a local $\ell$-translation of $\phi$ to
 {\sf GL}.  
\begin{itemize}
\item
Suppose $a$ is not in ${\sf bo}^\circ(\phi)$. We treat the case of a $\wedge$-occurrence, the other
cases being similar.
We take ${\sf js}^\star(\phi,a) := ({\sf js}^\star(\phi,S_0a)\wedge {\sf js}^\star(\phi,S_1a))$. 
\item
Suppose $a$ is in ${\sf bo}^\circ(\phi)$. 
 Then, ${\sf js}^\star(\phi,a) :=  {\sf js}^\star(\snip {(\sugr \phi a)} {\top},a)$.
 \item
 ${\sf js}(\phi) := {\sf js}^\star(\phi,r_\phi)$.
\end{itemize}

We have: 

\begin{lemma}
${\sf js}^\star(\phi,a)$ is in $\mathbb L$.
\end{lemma}

\begin{lemma}
${\sf js}^\star(\phi,a) = {\sf js}(\sugr\phi a)$.
\end{lemma}

\noindent
The desired results follow trivially by course of values induction on  $\mathsf c(\phi)$ and
by guard induction on ${\sf bo}^\circ(\phi)$.

%\begin{lemma}
%Suppose $\phi\in \mathbb L$. Then ${\sf js}(\phi) \bisim \phi$.
%\end{lemma}

We prove a result on commutation with substitution.

\begin{lemma}\label{subbismurf}
Let $\sigma$ be a substitution on $Q$.
Then, for $a\in V_\phi$, we have 
${\sf js}^\star (\phi \sigma, a) \bisim {\sf js}^\star(\phi,a)({\sf js}\circ \sigma)$.   
\end{lemma}

\begin{proof}
The proof is by course of values induction on  $\mathsf c(\phi)$ and by guard induction on ${\sf bo}^\circ(\phi)$.

Suppose $a\not\in {\sf bo}^\circ(\phi)$. Suppose, e.g., $a$ is a $\wedge$-occurrence.
We have:
\begin{eqnarray*}
{\sf js}^\star (\phi \sigma, a) & = & ({\sf js}^\star(\phi\sigma, S_0a) \wedge {\sf js}^\star(\phi\sigma,S_1a)) \\
& \bisim & ({\sf js}^\star(\phi,S_0a)({\sf js}\circ \sigma) \wedge {\sf js}^\star(\phi,S_1a)({\sf js}\circ \sigma)) \\
&\bisim & {\sf js}^\star(\phi,a)({\sf js}\circ \sigma)
\end{eqnarray*}
The cases where $a$ is an occurrence of a variable not in $Q$ or where $a$ is an occurrence of another connective
are similar. Suppose $a$ is an occurrence of $q\in Q$. We have:
\begin{eqnarray*}
 {\sf js}^\star (\phi \sigma, a) & \bisim & {\sf js}(\sugr {\phi\sigma} a) \\
 & = & {\sf js}({q\sigma}) \\
  & = & q ({\sf js}\circ \sigma) \\
 & = & {\sf js}(q) ({\sf js}\circ \sigma) \\
 & = & {\sf js}(\sugr \phi a) ({\sf js}\circ \sigma) \\
& \bisim &  {\sf js}^\star(\phi,a)({\sf js}\circ \sigma)
 \end{eqnarray*}
 
 Suppose $a\in {\sf bo}^\circ(\phi)$. We have:\qedright
 \begin{eqnarray*}
 {\sf js}^\star (\phi \sigma, a) & = &  {\sf js}^\star(\snip {(\sugr {\phi\sigma} a)} {\top},a) \\
 & = &  {\sf js}^\star(\snip {(\sugr \phi a)} {\top}\sigma,a) \\ 
 & \bisim & {\sf js}^\star(\snip {(\sugr \phi a)} {\top},a)({\sf js}\circ \sigma) \\
 & = &  {\sf js}^\star(\phi,a)({\sf js}\circ \sigma)
 \end{eqnarray*}
\end{proof}

We write $\equiv^\circ$ for provable equivalence in \glc\ and $\equiv$ for provable equivalence in {\sf GL}. 

%\begin{lemma}
%We have ${\sf js}(\phi) \equiv \phi$.
%\end{lemma}

%\begin{proof}
%We show by course of values induction on the number of cycles in $\phi$ and by guard induction on the
%$\opr$-occurrences of a cycle that ${\sf js}^\star(\phi,a) \equiv \sugr\phi a$.

%The case where $a$ is not a $\opr$-occurrence or $a$ is not on a cycle is easy. Suppose $a$ is an $\opr$-occurrence
%on a cycle. We have:
%\begin{eqnarray*}
 %{\sf js}^\star(\phi,a) & = &  {\sf js}^\star\snip {(\sugr \phi a)} \top \\
%& \equiv & \snip {(\sugr \phi a)} \top \\
%& \equiv & \sugr \phi a
%\end{eqnarray*}
%Here the second equivalence is by the induction hypothesis and the third 
%equivalence is by Lemma~\ref{kleutersmurf}.
%\end{proof}

For the proof of our main insight, Theorem~\ref{dicksmurf}, we need a lemma, that
is a strengthening of Theorem~\ref{unitrans}. Consider formulas $\phi$ and $\phi'$.
Let $F:V_\phi \to \mathbb L$ and let $F':V_{\phi'}\to \mathbb L$. 
We define ${\sf C}_\phi(F)$ as the conjunction of formulas $F(a) \iff \$(F(S_0a),\ldots,F(S_{n-1}a))$,
where $\$$ is the label of $a$ in $\phi$ and where $\$$ is $n$-ary. We treat the variable
as a 0-ary operation here.

\begin{lemma}\label{muzieksmurf}
Suppose $R$ is a bisimulation between $\phi$ and $\phi'$ and $aRa'$.
Then,
\begin{enumerate}[i.]
\item
$\dotbox {\sf C}_\phi(F) \wedge \dotbox{\sf C}_{\phi'}(F') \vdash_{\sf GL} F(a) \iff F(a')$.
\item
$\opr {\sf C}_\phi(F) \wedge \opr{\sf C}_{\phi'}(F') \vdash_{\sf GL} \opr(F(a) \iff F(a'))$.
\end{enumerate}
\end{lemma}

\begin{proof}
The proof of (i) is an immediate adaptation of the proof of Theorem~\ref{unitrans}, replacing
L\"ob's Rule by the Strengthened L\"ob's Rule. Item (ii) is immediate from (i).
\end{proof}

\begin{theorem}\label{dicksmurf}
${\sf js}^\star(\phi,\cdot)$ is a local translation of $\phi$ into {\sf GL}.
\end{theorem}

\begin{proof}
We employ course of values induction on the number of cycles in $\phi$ and, then,
guard induction on ${\sf bo}^\circ(\phi)$.
The only non-trivial case is where $a\in {\sf bo}^\circ(\phi)$. So, suppose $a\in {\sf bo}^\circ(\phi)$.
We have ${\sf js}^\star(\phi,a) = {\sf js}^\star (\snip {(\sugr\phi a)} \top, a)$.
Let us write $\psi_b := {\sf js}^\star(\phi,b)$, for $b\in V_\phi$  and $\psi'_c := {\sf js}^\star (\snip {(\sugr \phi a)} p, c)$, for 
$c \in V_{\snip {(\sugr \phi a)} p}$, where $p$ is a fresh variable. Let $S'$ be the successor function of $\snip {(\sugr \phi a)} p$.

We have (\dag):
\begin{eqnarray*}
 \psi_a & = &  {\sf js}^\star (\snip {(\sugr \phi a)} p [p:\top], a) \\
 & = &  \psi'_a [p:\top] \\
  & = & (\opr \psi'_{S'_0a}) [p:\top] \\
  & \equiv & (\opr \psi'_{S_0a}) [p:(\opr \psi'_{S'_0a}) [p:\top]] \\
  & = &  (\opr \psi'_{S'_0a}) [p:\psi_a] \\
  & = & \opr (\psi'_{S'_0a} [p:\psi_a]) 
  \end{eqnarray*}
  Here the second step is by Lemma~\ref{subbismurf}.
  The third step uses the fact that $a$ is not on a cycle in 
  $\snip {(\sugr \phi a)} p [p:\top]$.
  The fourth step uses Lemma~\ref{dejosmurf}.
  
  We now define 
  $F$ from $V_{\sugr \phi a}$ to $\mathbb L$ as follows.
  $F(a) := \psi_a$ and $F(b) := \psi'_b[p:\psi_a]$ if $a \neq b$.
  We prove that $F$ is a local translation from $\sugr \phi a$ to {\sf GL}.
  The case of $a$ is by (\dag). Suppose $b$ is a $\$$-occurrence unequal to $a$, where $\$$ is $n$-ary.
  By the induction hypothesis, we have ${\sf GL}\vdash \psi'_b \iff \$(\psi'_{S'_0b},\dots,\psi'_{S'_{n-1}b})$ and
  hence ${\sf GL} \vdash \psi'_b[p:\psi_a] \iff  \$(\psi'_{S'_0b}[p:\psi_a],\dots,\psi'_{S'_{n-1}b}[p:\psi_a])$.
  We note that:
  \begin{enumerate}[i.] 
  \item $ \psi'_b[p:\psi_a] = F(b)$;
  \item
  if $S'_ib$ is not a $p$-occurence, then   
  $\psi'_{S'_ib}[p:\psi_a] = \psi'_{S_ib}[p:\psi_a] = F(S_ib)$;
  \item
  if $S'_ib$ is a $p$-occurrence, then $\psi'_{S'_ib}[p:\psi_a] = \psi_a = F(a) = F(S_i b)$.
  \end{enumerate}
  
  \noindent
  Thus, we have ${\sf GL} \vdash {\sf C}_{\sugr\phi a}(F)$ and, hence, (\ddag) ${\sf GL} \vdash \opr {\sf C}_{\sugr\phi a}(F)$
  
  We define $G$ from   $V_\phi$ to the acyclic formulas by $G(b) := {\sf js}^\star(\phi,b)$.
  It follows that:
   \begin{eqnarray*}
   \opr {\sf C}_{\phi}(G)  & \vdash_{\sf GL} &    \opr {\sf C}_{\sugr \phi a}(F) \wedge \opr {\sf C}_{\phi}(G) \\
   & \vdash_{\sf GL} &  \opr  ( F(S_0a)\iff  G(S_0a))    \\
   & \vdash_{\sf GL} &  \opr F(S_0a) \iff \opr  G(S_0a) \\
   & \vdash_{\sf GL} & G(a) \iff \opr G(S_0a) 
   \end{eqnarray*}
   
   \noindent
   Here the first step is by (\ddag), the second step is by Lemma~\ref{muzieksmurf}(ii) and the fourth step is by (\dag). 
  
 Since $a$ was an arbitrary $\opr$-occurrence on a cycle, we find $\opr{\sf C}_\phi(G) \vdash_{\sf GL} {\sf C}_\phi(G)$. Ergo,
${\sf GL} \vdash {\sf C}(G)$. Thus, $G = {\sf js}^\star(\phi,\cdot)$ is a local translation of $\phi$ to {\sf GL}.
\end{proof}

\begin{theorem}\label{tuinsmurf}
{\sf js} carries an $\ell$-interpretation {\sf JS} of \glc\ in {\sf GL}.
\end{theorem}

\begin{proof}
By Theorem~\ref{dicksmurf}, ${\sf js}^\star(\phi,\cdot)$ is a local translation of $\phi$ in {\sf GL}.
It follows that {\sf js} is a translation of $\mathbb L^\circ$ into {\sf GL}. Our translation is an $\ell$-translation
by Lemma~\ref{subbismurf}. 

So we need just to verify the translations of the axioms and rules of \glc\ in {\sf GL}. However, the translations of
the axioms and rules of \glc\ except {\sf chl}5, are all instances of the same axioms and rules of {\sf GL}, modulo 
{\sf GL}-provable equivalence.
The axiom {\sf chl}5 follows by  Theorem~\ref{unitrans}.
\end{proof}

The identical translation {\sf emb} of $\mathbb L$ into \glc\ clearly carries an $\ell$-interpretation of {\sf GL} in \glc.
So, the pair {\sf JS}, {\sf Emb} forms an $\ell$-synonymy.
Since the arrows of a synonymy are faithful, it follows that \glc\ is conservative over {\sf GL}.
In other words, {\sf GL} is the acyclic fragment of \glc.       

Consider any acyclic $\phi$. Clearly, ${\sf GL} \vdash {\sf js}(\phi) \iff \phi$.
(Inspection shows that we even have ${\sf dj}(\phi) \bisim \phi$.)
Suppose $p$ is modalised in $\phi$.
We have: $\glc \vdash \digamma p.\phi \iff \phi[p:\digamma p.\phi]$.
It follows that
${\sf GL} \vdash {\sf js}( \digamma p.\phi ) \iff  \phi[p:{\sf js}(\digamma p.\phi)]$.
Thus ${\sf js}(\digamma p.\phi)$ is a de Jongh-Sambin explicit fixed point of 
$\phi$.

\begin{remark}
Of course, we could also develop the synonymy by using the known de Jongh-Sambin Theorem 
for {\sf GL}. However, the advantage of the present set-up is that we can see the algorithm
stated using guard recursion.
\end{remark}

\section{Further Inter-theory Relations}\label{intersection}
We develop the relations of {\sf CHL} to two other theories, to wit, Multiple Fixed Point Theory {\sf MFT} and
Well-foundedness Logic {\sf WfL} (or: ${\sf HL}(\mathbb L^\grullet)$).
In Subsection~\ref{mufithe}, we discuss {\sf MFT}. Our result on {\sf MFT} will play a role in the definition of
arithmetical interpretations. In Subsection~\ref{geaardesmurf}, we address {\sf WfL}.
Our result in that subsection is an ingredient of one possible Kripke completeness proof for
{\sf CHL}.

\subsection{Multiple Fixed Point Theory}\label{mufithe}
We define the following \emph{theory} {\sf MFT} in $\mathbb L$
extended with fresh constants.
We have {\sf K} plus L\"ob's Rule plus, for every system of equations $\mathcal E$ on a set of fixed-point 
variables $Q$ for the modal language
without constants, constants ${\sf c}_{\mathcal E,q}$ and axioms stating that these constants solve $\mathcal E$.
Let us say that an extension $\Theta$ of {\sf MFT} (in the same language)
is a \emph{strong extension} if it is closed under necessitation and L\"ob's Rule. 
Here we do not demand that $\Theta$ is closed under substitution: it is a theory not a logic.

We adapt the notion of local translation to the new setting in the following way. In the {\sf CHL}-{\sf MFT} direction
we allow formulas containing the constants. In the {\sf MFT}-{\sf CHL} connection we allow the constants to be interpreted
by formulas. Similarly, for global translations.

We have an immediate adaptation of Theorem~\ref{unitrans} to the slightly
modified setting that we
formulate here for completeness.

\begin{theorem}\label{uniekesmurf}
Suppose $R$ is a bisimulation between $\phi$ and $\phi'$ and suppose $\mathscr T$  is a local translation of $\phi$ in $\Theta$ and 
$\mathscr T'$ is a local interpretation of $\phi'$ in $\Theta$.  Suppose further that, for any variable $p$, if $a$ and $a'$ are occurrences of $p$, then
$\Theta \vdash \mathscr T(a) \iff \mathscr T'(a')$. Then, whenever $aRa'$, we have
  $\Theta \vdash \mathscr T(a) \iff \mathscr T'(a')$.
\end{theorem}

Our next order of business is to prove the existence of a local translation of $\phi$.
In Section~\ref{cyhelo}, we introduced the mapping {\sf E}. For convenience, we repeat the definition here.
Consider a formula $\phi$. We assign to each $\opr$-occurrence $a$ a propositional variable $q_a$, where the
$q_a$ are pairwise distinct and also distinct from the propositional variables of $\phi$. We map the
nodes of $\phi$ to formulas of the language of ordinary modal logic as follows:
\begin{itemize}
\item
${\sf E}_\phi(a) := q_a$, if $a$ is a $\opr$-occurrence.
\item
${\sf E}_\phi(a) := ({\sf E}_\phi(S_0a) \wedge {\sf E}_\phi(S_1a))$, if $a$ is a $\wedge$-occurrence.
Similarly, for the other connectives and for the propositional variable-occurrences in $\phi$.
\end{itemize}

We write $\psi_a$ for ${\sf E}_\phi(a)$. We note that $\psi_a$ is cycle-free.

Let $\mathcal E$ be defined by $q_a\mathcal E := \opr {\sf E}_\phi(S_0a)$, where $a$ is an $\opr$-occurrence in $\phi$. 
In {\sf MFT} this system of equations has
a solution, say $q_a\mathcal F := c_a$.
Finally we define, for every $b \in V_\phi$, the mapping ${\sf cyco}^\star_\phi$ by ${\sf cyco}^\star_\phi(b) : \psi_b\mathcal F$.

\begin{lemma}
${\sf cyco}^\star_\phi$ is a local translation of $\phi$ in {\sf MFT}. 
\end{lemma}

\begin{proof}
The cases where $b$ is not a $\opr$-occurrence are simple. 
Suppose $b$ is an $\opr$-occurrence. We need that
${\sf MFT} \vdash \psi_b\mathcal F \iff \opr \psi_{S_0b}\mathcal F$.
However, this is precisely ${\sf MFT} \vdash c_b \iff \opr \psi_{S_0b}\mathcal F$,
the promised solution of $\mathcal E$.
\end{proof}

It follows that {\sf cyco} defined by ${\sf cyco}(\phi) := {\sf cyco}^\star_\phi(r_\phi)$ is a global
translation of $\mathbb L^\circ$ in {\sf MFT}. Using Theorem~\ref{uniekesmurf}, we now find:

\begin{theorem}
There is an interpretation {\sf CyCo} based on {\sf cyco} of {\sf CHL} in {\sf MFT}.
\end{theorem}

In the other direction we define a translation {\sf cocy} that commutes with propositional variables 
and connectives and that sends a constant $c_q$ in introduced for a system of equations $\mathcal E$
to $\digamma_q\mathcal E$, the solution of $\mathcal E$ for $q$ as guaranteed by Theorem~\ref{mufipohe}.
We find:

\begin{theorem}
There is an interpretation {\sf CoCy} based on {\sf cocy} of {\sf MFT} in {\sf CHL}. Moreover, this interpretation is
unique.
\end{theorem}

\begin{proof}
The verification that {\sf cocy} carries an interpretation is entirely as expected. For the uniqueness we use
Theorem~\ref{gejosam}. 
\end{proof}

Using the analogue of Theorem~\ref{gejosam} in {\sf MFT}, we find:

\begin{theorem}
{\sf CyCo} and {\sf CoCy} form a synonymy.
\end{theorem}

We note that the synonymy we are looking at here is a synonymy of \emph{theories} not logics.
It is a form of sameness weaker than the sameness of ${\sf GL}^\circ$ and {\sf GL}.

\subsection{Well-Foundedness Logic}\label{geaardesmurf}
It is interesting to note that the global translations {\sf JS} and {\sf Emb} still carry a synonymy between {\sf WfL} plus
{\sf L}3 and  {\sf GL} when we restrict {\sf JS} to $\mathbb L^\grullet$. In this case we simply have
${\sf JS}(\opr^\grullet \phi) = \opr({\sf JS}(\phi) \wedge \top)$.

We borrow a result form a forthcoming paper with Tadeusz Litak:
we have the completeness theorem for {\sf WfL} in finite acyclic Kripke models.

We remind the reader that ${\sf bo}(\phi)$ is the set of $\opr$-occurrences in $\phi$.
Modulo the announced result, we have:

\begin{antheorem}
Consider any $\mathbb L$-formula $\phi$. We use $\psi_a :={\sf E}_\phi(a)$ as in Subsection~\ref{mufithe}.
 We have:
\[ {\sf CHL} \vdash \phi \text{ iff }\bigwedge_{a\in {\sf bo}(\phi)} \dotbox^\grullet(q_a\iff \opr \psi_{S_0a}) \vdash_{\sf WfL} \psi_r.\]
\end{antheorem} 

\begin{proof}
Suppose $\bigwedge_{a\in {\sf bo}(\phi)} \dotbox^\grullet(q_a\iff \opr \psi_{S_0a}) \vdash_{\sf WfL} \psi_r$. Then, \emph{a fortiory},
\[\bigwedge_{a\in {\sf bo}(\phi)} \dotbox^\grullet(q_a\iff \opr \psi_{S_0a}) \vdash_{\sf CHL} \psi_r.\]
Then, by Theorem~\ref{thed}, we have: (\dag)
$\bigwedge_{a\in {\sf bo}(\phi)} \dotbox^\grullet(q_a\iff \opr \psi_{S_0a}) \vdash_{\sf CHL} \phi$.
We note that $\phi$ does not contain occurrences of the $q_a$. By the Simultaneous Fixed Point Theorem~\ref{mufipohe},
we have {\sf CHL}-verifiable solutions $\chi_a$ of the equations in the antecedent conjunction. 
Substituting these for the $q_i$, we find ${\sf CHL}\vdash \phi$.

Conversely, suppose $\bigwedge_{a\in {\sf bo}(\phi)} \dotbox^\grullet(q_a\iff \opr \psi_{S_0a}) \nvdash_{\sf WfL} \psi_r$. 
By the announced Completeness Theorem for {\sf WfL}, we have an acyclic finite Kripke model $\mathcal K$ with
root $\mathfrak r$ such that $\mathfrak r\Vdash \bigwedge_{a\in {\sf bo}(\phi)} \dotbox^\grullet(q_a\iff \opr \psi_{S_0a})$
and $\mathfrak r \nVdash \psi_r$. Clearly the forcing relation for
$\mathbb L^\grullet$ is the restriction of the forcing relation on $\mathbb L^\circ$ to $\mathbb L^\grullet$.
Since $\mathcal K$ forces {\sf CHL} it follows, by Theorem~\ref{thed}, that $\mathfrak r \nVdash \phi$.
So, ${\sf CHL} \nvdash \phi$.
\end{proof}

\begin{antheorem}
We have completeness for {\sf CHL} in finite acyclic Kripke models.
\end{antheorem} 

\begin{proof}
Suppose ${\sf CHL} \nvdash \phi$. Then,
$\bigwedge_{a\in {\sf bo}(\phi)} \dotbox^\grullet(q_a\iff \opr \psi_{S_0a}) \nvdash_{\sf WfL} \psi_r$.
So, we have, by the announced completeness theorem for {\sf WfL}, a Kripke model $\mathcal K$ with root $\mathfrak r$ such that 
$\mathfrak r\Vdash \bigwedge_{a\in {\sf bo}(\phi)} \dotbox^\grullet(q_a\iff \opr \psi_{S_0a})$
and $\mathfrak r \nVdash \psi_r$. It follows by Theorem~\ref{thed}, that $\mathfrak r \nVdash \phi$.
\end{proof}

\section{Arithmetical interpretations}\label{arithintsection}
In this section, we introduce the notion of \emph{arithmetical interpretation} and verify its basic properties.

We fix a theory $U$, an interpretation $N$ of ${\sf S}^1_2$ in $U$. We suppose that we have a provability predicate
{\sf bew} for numerals in $N$ that satisfies {\sf L}1,2 plus L\"ob's Rule. We write $\graysq A$ for ${\sf bew}(\gnum A)$. 

\begin{remark}
In our paper \cite{viss:look19}, we show that if we have a predicate 
that uniformly semi-represents a given axiom set of $U$ (w.r.t. $N$) in $U$.
The predicate ${\sf prov}^N_{[\alpha]}$, where {\sf prov} is a standard provability predicate, has
the desired properties we ask of {\sf bew}. 
We discuss these ideas in some detail in Appendix~\ref{aritexa}.
\end{remark}

%We treat the syntactic operations like {\sf conj} as terms. 
%These are to be eliminated in the usual way.

A \emph{translation} of $\phi$ in $U$ (for {\sf bew}) is a mapping $\mathscr T$ from $V_\phi$ to $U$-sentences, such that:
\begin{itemize}
\item
If $a$ is an occurrence of $\wedge$ then $U \vdash \mathscr T(a) \iff (\mathscr T(S_0a) \wedge \mathscr T(S_1a))$.
Similarly for the other non-box connectives. 
\item
If $a$ is an occurrence of $\opr$ and $Sa = \tupel{b}$, then $U \vdash \mathscr T(a) \iff \graysq \mathscr T(S_0a)$.
\item
If $a$ and $a'$ are occurrences of $p$, then $U \vdash \mathscr T(a) \iff \mathscr T(a')$.
\end{itemize}
An \emph{interpretation} of $\mathbb L^\circ$ in $U$ is a mapping $\mathscr U$ from  $\mathbb L^\circ$
to $U$-sentences that commutes modulo $U$-provability with the propositional connectives and the commutes 
with $\graysq$ in the sense that $U \vdash \mathscr U(\opr\phi) \iff \graysq \mathscr U (\phi)$.

Modulo $U$-provable equivalence, arithmetical translations are preserved under bisimulation.

\begin{theorem}\label{uniekesmurftwee}
Suppose $R$ is a bisimulation between $\phi$ and $\phi'$ and suppose $\mathscr T$  is a translation  of $\phi$ in $U$
and $\mathscr T'$ is a translation of $\phi'$ in $U$.  Suppose further that, for any variable $p$, if $a$ and $a'$ are occurrences of $p$, then
$U \vdash \mathscr T(a) \iff \mathscr T'(a')$. Then, whenever $aRa'$, we have
  $U \vdash \mathscr T(a) \iff \mathscr T'(a')$.
\end{theorem}

The proof is just a minor variation of the proof of Theorem~\ref{unitrans}. We note that it follows that,
if we have translations $\mathscr T_\phi$ for each $\phi\in \mathbb L^\circ$, we can base an interpretation
on them by defining $\mathscr U(\phi) := \mathscr T_\phi(r_\phi)$. The Uniqueness Theorem
guarantees that the local pieces add up to a coherent whole.

As a preparation of the theorem concerning the existence of arithmetical translations, we
first remind the reader of the simultaneous fixed point lemma and provide two proofs of it.

\begin{lemma}
Suppose $N$ interprets ${\sf S}^1_2$ in $U$.\footnote{We assume 
here that the variables of a formula or term are among the variables that are displayed.} 
Consider $U$-formulas $A_i(x_0,\dots x_{n-1},\vec y\,)$ for $i<n$.
Here the variables $x_i$ range over the domain of $N$.
 Then, there are formulas $B_0(\vec y\,),\dots,B_{n-1}(\vec y\,)$ such that
\[U\vdash B_i(\vec y\,) \iff A_i(\gnum{B_0(\dot{\vec y}\,)},\dots,\gnum{B_{n-1}(\dot{\vec y}\,)},\vec y\,),\] for $i<n$.\footnote{We employ
the Smory\'nsky dot notation.}
\end{lemma}

\begin{proof}[First Proof.]
The proof is by induction on $n$. By the usual fixed point lemma with parameters we find a
$C_0(x_1,\dots, x_{n-1},\vec y\,)$ such that 
\[ U \vdash C_0(x_1,\dots, x_{n-1},\vec y\,) \iff A_0(\gnum{C_0(\dot x_1,\dots, \dot x_{n-1},\dot{\vec y}\,)},x_1,\dots,x_{n-1},\vec y\,).\] 
Now we define:
\[ C_{i+1}(x_1,\dots, x_{n-1},\vec y\,):=  A_{i+1}(\gnum{C_0(\dot x_1,\dots, \dot x_{n-1},\dot{\vec y}\,)},x_1,\dots,x_{n-1},\vec y\,).\] 
We solve the system $C_1,\dots C_{n-1}$ as is guaranteed by the induction hypothesis, resulting in $B_1,\dots,B_{n-1}$.
Finally we set:
\[B_0(\vec y\,) := C_0(\gnum{B_1(\dot{\vec y}\,)},\dots, \gnum{B_{n-1}(\dot{\vec y}\,)},\vec y\,).\]
It is now easy to see that the $B_i$ form the desired simultaneous fixed point.
\end{proof}

\begin{proof}[Second Proof.]
Let us write $x[w]$ for the term representing the result of substituting the numeral of $w$ for the variable $z$ in the formula represented by $x$.
We define $D(z,x,\vec y\,) := \bigvee_{i<n}(z= \underline i \wedge A_i(x[\underline 0],\dots, x[\underline{n{-}1}],\vec y\,))$.
Let $E(z,\vec y\,)$ be the fixed point of $D$ w.r.t. $x$. Then, we can take
$B_i(\vec y\,) := E(\underline i,\vec y\,)$.
\end{proof}

The notion of interpretation for {\sf MFT} in $U$ is the obvious one, since {\sf MFT} has acyclic syntax.
Using the simultaneous fixed point lemma, we find:
\begin{theorem}
For every $f$ from the propositional variables $U$-sentences,
there is an interpretation ${\sf FA}_f$ based on ${\sf fa}_f$ of {\sf MFT} in  $U$ such that
${\sf fa}_f(p) = f(p)$. 
\end{theorem}

\begin{theorem}\label{existentiesmurf}
For every $f$ from the propositional variables $U$-sentences,
there is an interpretation ${\sf HA}_f$ based on a translation ${\sf ha}_f$ of {\sf CHL} in  $U$ such that
${\sf fa}_f(p) = f(p)$. This interpretation  is unique modulo provable equivalence.
\end{theorem}

\begin{proof}
We take ${\sf ha}_f := {\sf fa}_f \circ {\sf cyco}$. For the verification of {\sf chl}5 and for  uniqueness,
we use Theorem~\ref{uniekesmurftwee}.
\end{proof}

We note that we can view ${\sf ha}_f$ as a mapping of formulas modulo bisimulation to
elements of the diagonalised algebra of $U$.

If $U$ satisfies some further schematic logical principle, it is easy to see that ${\sf ha}_f$ carries
an interpretation $H$ of {\sf CHL} plus that logical principle. So, e.g., if $\graysq$ satisfies {\sf L}3 in $U$,
we have in interpretation $H$ based on ${\sf ha}_f$ of ${\sf GL}^\circ$ in $U$, etcetera.

\begin{remark}
Clearly, we can use Theorem~\ref{tuinsmurf}  to prove arithmetical completeness for ${\sf GL}^\circ$
 in a $\Sigma^0_1$-sound extension $U$ of Elementary Arithmetic for Fefermanian provability with respect
to an elementary representation $\alpha$ of the axiom set. We can do it directly.
It is somewhat remarkable that (a lifted version of) Theorem~\ref{existentiesmurf} is not needed in the proof.
The Solovay construction delivers  translations with the desired properties directly. 

Here is a sketch of how this works.
Consider any tail-model $\mathcal K$ for ${\sf GL}/{\sf GL}^\circ$. See \cite{viss:prov84} for this notion. 
Let $X$ be a finite or cofinite set of nodes. We write $[X]$ for
(an appropriate paraphrase of) $\exists x \in X \, \ell = x$, where $\ell$ is the limit statement constructed
by Solovay (for $\mathcal K$ and for $\opr_\alpha$). We can now show:
$U \vdash [X\cap Y] \iff [X] \wedge [Y]$ and, similarly, for the other non-modal connectives.
We write \[\opr X := \verz{y\mid \forall x \, (y\sqsubset x \To x\in X)}.\] We have:
$U \vdash [\,\opr X] \iff  \opr_\alpha [X]$. 

Consider any $\phi\in \mathbb L^\circ$.
We define  ${\sf S}^\star_{\mathcal K,\phi}(a) := [\verz{x \mid x\Vdash \sugr \phi a}]$.
(Here we need to check that $\verz{x \mid x\Vdash \psi}$ is always finite or co-finite not just for 
$\mathbb L$ as is proven in  \cite{viss:prov84}, but also for $\mathbb L^\circ$.)
It is easy to see that we can base an interpretations ${\sf S}_{\mathcal K}$ on the 
${\sf S}^\star_{\mathcal K,\phi}$ and that these interpretations witness the desired arithmetical
completeness for ${\sf GL}^\circ$.

Thus, in this proof we use the Magari Algebra of $\mathcal K$ to replace the use of {\sf MFT} in the proof of
Theorem~\ref{existentiesmurf}. We interpret ${\sf GL}^\circ$ in the Magari Algebra of $\mathcal K$ via
our result on evaluation in Kripke models and we interpret $\mathcal K$ in $U$ via the standard
Solovay argument as applied to tail models.
\end{remark}

%\bibliographystyle{alpha}
%\bibliography{provint}

\appendix

\section{Arithmetical Matters}\label{arima}
In this appendix, we summarise the relevant ideas from our paper \cite{viss:look19} and
prove some additional results that connect that paper to the present one.

In Subsection~\ref{gargamel}, we provide some basic insights. In Subsection~\ref{alcibiadessmurf},
we reprove the version of the Second Incompleteness Theorem from \cite{viss:look19} using the notations
of the present paper. In Subsection~\ref{tracl}, we revisit the transitive closure modality $\opr^\grullet$ and
show that some of its properties also hold globally. Moreover, we show that, under certain conditions,
the Fefermanianness of the input modality is preserved to its transitive closure. 
In Subsection~\ref{smoothsmurf}, we discuss the operations of Craigification and Smoothening.
Finally, we consider an example with some remarkable properties in Subsection~\ref{aritexa}.

\subsection{Preliminaries}\label{gargamel}
 We consider theories $T$ in predicate logic of finite signature.
We allow $T$ to have any complexity. We always assume that $T$ is equipped with
an interpretation $N:{\sf S}^1_2\to T$. Unless stated otherwise, displayed variables 
 range over the domain of $N$. In other words, we pretend that $T$ is an arithmetical
theory. Moreover, we will use single variables to range over the domain of $N$, even if $N$ might be
multi-dimensional.

Our treatment is, at places, somewhat dependent on details of the coding.
 Let us simply assume that we base our coding on an alphabet that contains (at least)
 the signs of the language plus some extra brackets ``$[$'' and ``$]$''.
 Our  G\"odel numbering is given by the length-first ordering of strings from this alphabet.
 We code finite sets of expressions as (the G\"odel numbers of) strings of expressions, so $\gnum{A_0A_1\dots A_{k-1}}$
 and finite sets of numbers as the finite sets of their efficient numerals.
 We represent proofs with assumptions as  strings of formulas where some formulas, the assumptions, are
 enclosed between square brackets.  
 We code finite functions $B_i \mapsto n_i$ from sentences to numbers as strings of 
the form $B_0\overline n_0 \dots B_{m-1}\overline n_{m-1}$. Etcetera.
Of course, most of these details are immaterial. The main things we need are
properties like the following: the function that sends $n$ to the G\"odel number of its
efficient numeral, $\gnum{\overline n}$, is p-time; the assumption set of a proof can be efficiently extracted from
the proof; etcetera.

Consider a $T$-predicate {\sf bew}. 
We write $\graysq A$ for ${\sf bew}(\gnum A)$. Here the numerals are defined w.r.t.\ $N$.
We say that {\sf bew} is a \emph{provability predicate} if it satisfies the first two L\"ob conditions
w.r.t. $T$.
\begin{enumerate}[{\sf L}1.]
\item
If $T \vdash A$, then $T \vdash \graysq A$.
\item
$T \vdash (\graysq A \wedge \graysq\, (A \to B)) \to \graysq B$.
\end{enumerate} 

 A special class of provability predicates are the Feferman-style predicates. This works as follows.
 We fix  arithmetisations of ${\sf proof}^\dag(p,y)$ and ${\sf ass}(p)$, where  ${\sf proof}^\dag(p,y)$ is a good arithmetisation
 of:
$p$ codes a proof in predicate logic of $y$ (of the ambient signature) from assumptions in ${\sf ass}(p)$.
We write:
\begin{itemize}
\item
 ${\sf proof}_\alpha(p,y)$ for $({\sf proof}^\dag(p,y) \wedge \forall z\in {\sf ass}(p)\, \alpha(z))$, 
 \item
${\sf prov}_\alpha(y)$ for $\exists p\,{\sf proof}_\alpha(p,y)$,
\item
$\opr_\alpha A$ for ${\sf prov}_\alpha(\gnum{A})$.
\end{itemize}
Now consider the theory $T$ with axiom set $X$. We suppose we have an interpretation $N$ of ${\sf S}^1_2$ in $T$.
Suppose that $\alpha$ semi-numerates $X$ in $T$, i.e., if $A\in X$, then $T \vdash \alpha(\gnum A)$.
In these circumstances, it is easy to see that ${\sf prov}_\alpha^N$ is a provability predicate for $T$.
We say that ${\sf prov}_\alpha^N$ is a Fefermanian provability predicate. 
We note that $\alpha$ need not be of the form $\beta^N$, where $\beta$ is an arithmetical
predicate.

If $X_0$ is a finite set of $T$-sentences, we write $[X_0]$ for 
$\bigvee_{B\in X_0} x=\gnum{B}$. We write $[A]$ for $[\verz{A}]$.
We note that $\opr_{[\emptyset]}$, $\opr_{[\top]}$ and $\opr_\bot$ all represent provability in predicate logic
of the ambient signature.

%In the definition of ${\sf L}^{\sf glob}1$ we may, equivalently, 
%replace ${\sf prov}({\sf imp}(\gnum A,b))$, by the more
%readable ${\sf prov}_{[A]}(b)$. 

Here are some further definitions.
\begin{itemize}
\item
${\sf bew}$ is an \emph{LR-provability predicate} for $T$ if, it satisfies ({\sf L}2) and $T$ is closed under
L\"ob's rule {\sf LR} for {\sf bew}: if $T \vdash \graysq A \to A$, then $T \vdash A$.
\item
${\sf bew}$ is a \emph{uniform provability predicate} for $T$, if it satisfies the following three principles
\begin{enumerate}[${\sf L}^{\sf un}$1.]
\item
Whenever $T \vdash A$, there is finite set of $T$-sentences $X_0$ such that, $X_0\vdash A$ and, for each $B\in X_0$, we have
 $T \vdash B$ and $X_0 \vdash \graysq B$.
 \item
 $T \vdash \forall b \, ({\sf prov}^N_{[\emptyset]}(b) \to {\sf bew}(b))$.
  \item
$T\vdash  \forall a \,\forall b\, (({\sf bew}(a) \wedge {\sf bew}({\sf imp}(a,b))) \to {\sf bew}(b))$.
 \end{enumerate}
 \item
 {\sf bew} is a \emph{global provability predicate}, if it satisfies {\sf L1} and ${\sf L}^{\sf un}2$ and  ${\sf L}^{\sf un}3$.
 \end{itemize}
 
 We provide some basic results on these notions.
 
 \begin{theorem}
 Suppose {\sf bew} is an LR-provability predicate for $T$. Then, {\sf bew} is a provability predicate for $T$.
 \end{theorem}
 
 \begin{proof}
 The proof is just a minor adaptation of an argument due to Dick de Jongh. Suppose $\graysq$ is an LR-provability predicate for $T$.
 Suppose $T \vdash A$. Then, by {\sf K}-reasoning, $T \vdash \graysq (A \wedge \graysq A) \to (A \wedge \graysq A)$.
 So, by L\"ob's rule, we have $T \vdash \graysq A$.
 \end{proof}
 
 \begin{theorem}
If {\sf bew} is a uniform provability predicate, it is a global provability predicate. If {\sf bew} is a global provability predicate, it is
a provability predicate.
 \end{theorem}
 
 \begin{proof}
 We verify e.g. that a uniform provability predicate satisfies {\sf L1}. Suppose
 $T \vdash A$. We find $X_0 = \verz{B_0,\dots, B_{n-1}}$ as promised. We have $X_0 \vdash A$ and, hence,
 $T \vdash \opr^N_{[\emptyset]}C$, where $C :=(B_0 \to (B_1 \to \dots  (B_{n-1}  \to A)\dots ))$.  
 It follows that $T\vdash \graysq C$. We also have $T \vdash \graysq B_i$ for $i<n$, so by repeated
 application of ${\sf L}^{\sf un}3$, we find $T \vdash \graysq A$.
 \end{proof}
 
 The next theorem is obvious.
 
\begin{theorem}
Any Fefermanian provability predicate is global.
\end{theorem}  

\subsection{A Version of the Second Incompleteness Theorem}\label{alcibiadessmurf}
We present the relevant result of  \cite{viss:look19} in the terminology of the present paper. 

\begin{theorem}\label{grotesmurf}
Suppose $\graysq$ is a uniform provability predicate for $T$.
Then, $\graysq$ is an LR-provability predicate for $T$.
\end{theorem}

\begin{proof}
 Let $\graysq$ be a uniform provability predicate for $T$.
 Suppose $T \vdash \graysq A \to A$. Let $C$ be the conjunction of
the following statements:
\begin{itemize}
\item 
$(\bigwedge{\sf S}^1_2)^N$ (we assume that the axioms of ${\sf S}^1_2$ include the axioms of identity),
\item  
 $\forall b \, ({\sf prov}^N_{[\emptyset]}(b) \to {\sf bew}(b))$,
 \item
 $ \forall a \,\forall b\, (({\sf bew}(a) \wedge {\sf bew}({\sf imp}(a,b))) \to {\sf bew}(b))$,
 \item
 $\graysq A \to A$.
 \end{itemize}
 Let $X_0$ be as promised for $C$ in the definition of uniformity.
We have  $X_0\vdash  C$.
Hence,
\begin{eqnarray*}
 X_0\;\;  \vdash  \opr^N_{[X_0]}A &\to& \opr_{[\emptyset]}^N (\,\bigwedge X_0 \to A) \\
 & \to &  \graysq\, (\bigwedge X_0 \to A) \\
 & \to & \graysq\, A \\
 & \to & A
 \end{eqnarray*}
 So, we find that $X_0 \vdash \opr^N_{[X_0]}A \to A$. Since we have L\"ob's Rule for finitely axiomatised
 theories with standard axiomatisation (w.r.t. an interpretation $N$ of ${\sf S}^1_2$), 
 we find $X_0 \vdash A$ and, hence, $T \vdash A$.
\end{proof}

\begin{question}
It would be interesting to have an example of a provability predicate for which we have
the Second Incompleteness Theorem, but not L\"ob's Rule. 

It would be interesting to have an example of an LR-provability predicate  that is not
uniform.
\end{question}

\subsection{Global Properties of the Transitive Closure Modality}\label{tracl}
We show that a number of desirable properties can be lifted from {\sf bew} to
${\sf bew}^\grullet$.

Suppose {\sf bew} is an LR-provability predicate for $T$ w.r.t. $N$, where $N$ is an interpretation of
${\sf S}^1_2$. Let $a,b,c,\dots$ range over codes of sentences. We use {\sf conj}, {\sf imp}, etcetera,
for the arithmetisations of the obvious syntactical operations.

By the fixed point lemma with parameters we find ${\sf bew}^\grullet$, such that:
\begin{itemize}
\item
$T \vdash\; {\sf bew}^\grullet(a) \;\;\iff \;\; {\sf bew}({\sf conj}(a, \gnum{{\sf bew}^\grullet(\dot a)}))$.
\end{itemize}

\noindent
This is just our previous definition but now not locally for each sentence but for all sentences at once.
The definition is unproblematic using the fixed point construction with parameters.
We note that by the uniqueness result our global definition will locally coincide with whatever way
we implemented the local definitions.

\begin{theorem}\label{globalsmurf}
Suppose {\sf bew} is a global LR-provability predicate for $T$ w.r.t. $N$, where $N$ is an interpretation of
${\sf S}^1_2$.  Then,
${\sf bew}^\grullet$ is a global LR-provability predicate w.r.t. $N$.
\end{theorem} 

\begin{proof}
We assume the conditions of the theorem.
We can copy the reasoning of the proof of Theorem~\ref{stevigesmurf}, to show that
$\graysq^\grullet$ is closed under necessitation and L\"ob's Rule.

We verify ${\sf L}^{\sf un}2$ for ${\sf bew}^\grullet$.
 Let $A := (\bigwedge{\sf S}^1_2)^N$, where we assume ${\sf S}^1_2$ to include the
theory of identity. Reason in $T$. Suppose ${\sf prov}_{[\emptyset]}^N(b)$. By the formalisation in
${\sf S}^1_2$ of global $\exists \Sigma_1^{\sf b}$-completeness, we find 
${\sf prov}_{[\emptyset]}({\sf imp}(\gnum A,\gnum{{\sf prov}_{[\emptyset]}^N(\dot b)}))$.
It follows that ${\sf bew}({\sf imp}(\gnum A,\gnum{{\sf prov}_{[\emptyset]}^N(\dot b)}))$. Since, by necessitation, we also have
${\sf bew}(\gnum{A})$. So, we find, by ${\sf L}^{\sf un}3$ for {\sf bew}, that
 ${\sf bew}(\gnum{{\sf prov}_{[\emptyset]}^N(\dot b)})$. We leave $T$.
 We have shown:
 \[ T \vdash \forall b \, ({\sf prov}_{[\emptyset]}^N(b) \to {\sf bew}(\gnum{{\sf prov}_{[\emptyset]}^N(\dot b)})).\]
It follows that:
\begin{eqnarray*}
\graysq\,(\forall a\,({\sf prov}_{[\emptyset]}^N(a) \to {\sf bew}^\grullet(a))) \vdash_T \;\; {\sf prov}_{[\emptyset]}^N(b)
&  \to & ({\sf bew}(b) \wedge {\sf bew}(\gnum{{\sf prov_{[\emptyset]}}^N(\dot b)})) \\
& \to & ({\sf bew}(b) \wedge {\sf bew}(\gnum{{\sf bew}^\grullet(\dot b)})) \\
& \to & {\sf bew}({\sf conj}(b,\gnum{{\sf bew}^\grullet(\dot b)})) \\
& \to & {\sf bew}^\grullet(b).
\end{eqnarray*}
We may conclude that $\graysq(\forall a\,({\sf prov}_{[\emptyset]}^N(a) \to 
{\sf bew}^\grullet(a))) \vdash_T \forall a\,({\sf prov}_{[\emptyset]}^N(a) \to {\sf bew}^\grullet(a))$, hence, by L\"ob's Rule,
$T \vdash \forall a\,({\sf prov}_{[\emptyset]}^N(a) \to {\sf bew}^\grullet(a))$.

We verify ${\sf L}^{\sf un}3$. Let
$B := \forall a \,\forall b\, (({\sf bew}^\grullet(a) \wedge {\sf bew}^\grullet({\sf imp}(a,b))) \to {\sf bew}^\grullet(b))$.
We have:
{\footnotesize
\begin{eqnarray*}
\graysq B \vdash ({\sf bew}^\grullet(a) \wedge {\sf bew}^\grullet({\sf imp}(a,b))) & \to & 
{\sf bew}({\sf conj}(a,\gnum{{\sf bew}^\grullet (\dot a)})) \;\wedge\\
&&
{\sf bew}({\sf conj}({\sf imp}(a,b),\gnum{{\sf bew}^\grullet ({\sf imp}(\dot a,\dot b))}))  \\  
& \to & {\sf bew}(a) \wedge {\sf bew}({\sf imp}(a,b)) \; \wedge\\
&& {\sf bew}( {\sf conj}({\sf bew}^\grullet(\dot a), {\sf bew}^\grullet({\sf imp}(\dot a,\dot b) )))\\ 
&\to & {\sf bew}(b) \wedge {\sf bew}({\sf bew}^\grullet(\dot b)) \\
& \to &  {\sf bew}({\sf conj}(b,{\sf bew}^\grullet(\dot b))) \\
& \to & {\sf bew}^\grullet(b)
\end{eqnarray*}
}
It follows that $\graysq B \vdash_T B$.
So, by L\"ob's Rule, we have $T \vdash B$.
\end{proof}

We can now show that the $(\cdot)^\grullet$-operation preserves Fefermanianness.

\begin{theorem}\label{preserveersmurf}
Suppose ${\sf prov}_\alpha^N$ is an LR provability predicate. Then $({\sf prov}_\alpha^N)^\grullet$ is also a Fefermanian
provability predicate.
\end{theorem}

\begin{proof}
To simplify the presentation, we will omit the relativisation to $N$. We reason inside $N$ but for the fact that $\alpha$ is
not an internal $N$-formula.
Inspecting the fixed point construction, we see that
${\sf prov}_\alpha^\grullet(a)$ is of the form (\$) $\exists b\, (S_0(a,b) \wedge {\sf prov}_\alpha({\sf conj}(a,b)))$, where $S_0$ 
is an $\exists\Sigma^{\sf b}_1$-formula that represents the relevant term, such that, $T$-verifiably,
$\forall a\, \forall b\,(S_0(a,b) \iff b=  \gnum{{\sf prov}_\alpha^\grullet(\dot a)})$.
We will treat the existential quantifier of ${\sf prov}_\alpha$ in (\$) as giving the primary witness of
$ {\sf prov}_\alpha^\grullet(a)$. 

To define  $\alpha^\grullet(a)$, we, sloppily, use meta-notations in the object language.
Thus we write, for example, $(A\wedge B)$ for ${\sf conj}(a,b)$. I think the gain in readability outweighs
the loss of precision. The reader just should remember that, locally, the roman capitals represent internal
variables.

We take $\alpha^\grullet(A)$ iff $A$ is of the form $({\overline p=\overline p\,} \wedge B)$, where overlining gives
us efficient numerals and where $p$ is a primary witness of $\opr_\alpha^\grullet B$. In other words,
$p$ is a witness of $\opr_\alpha(B \wedge \opr^\grullet_\alpha B)$. The attentive reader will see that
the definition of $\alpha^\grullet$ is an incarnation of Craig's trick.

We work in $T$.
Suppose $\opr_\alpha^\grullet B$. Let $p$ be a primary witness. Then, $(\overline p=\overline p \wedge B)$ is in $\alpha^\grullet$.
Hence trivially $\opr_{\alpha^\grullet}B$.

Now suppose $\opr_{\alpha^\grullet}B$. Let $q$ be a witnessing proof.
Suppose $(\overline p_0 = \overline p_0 \wedge C_0)$, \dots,
$(\overline p_{k-1} = \overline p_{k-1} \wedge C_{k-1})$ are the
(possibly non-standardly many) axioms used in $q$.
We want to prove $\opr_\alpha^\grullet B$.
Let  $D := (C_0 \to (C_1 \to \dots (C_{k-1} \to B)\dots))$.
We claim that, inside  $\opr_\alpha$, we have:
\begin{enumerate}[I.]
\item
B
\item
$\opr_\alpha^\grullet C_i$, for $i<k$,
\item
 $\opr_\alpha^\grullet D$,
 \item
 $\forall a,b\, (({\sf prov}^\grullet_\alpha(a)  \wedge {\sf prov}^\grullet_\alpha({\sf imp}(a,b))) \to {\sf prov}_\alpha^\grullet(b))$.
 \end{enumerate}
 (I) follows from the fact that we can effectively transform $q$ into an $\alpha$-proof of $B$.
  We have (II) by the definition of the $p_i$. We note that the finite set of witnesses of the $\opr_\alpha\opr_\alpha^\grullet C_i$  is not much larger than $q$.
 (III)  follows from the fact that we can transform $q$ effectively into a witness of 
 $\opr_{[\emptyset]} D$. Since, by Theorem~\ref{globalsmurf}, ${\sf prov}_\alpha^\grullet$ is global,
it follows that $\opr_\alpha^\grullet D$ and, hence, $\opr_\alpha \opr_\alpha^\grullet D$. 
Finally, we have (IV) since $\opr_\alpha$ is global combined with necessitation for $\opr_\alpha$.

Now we combine (II), (III) and (IV) to effectively find a witness of $\opr_\alpha\opr_\alpha^\grullet B$.
So, combining this with (I), we find $\opr_\alpha (B\wedge \opr_\alpha^\grullet B)$.
Hence, we have $\opr_\alpha^\grullet B$. 
 \end{proof}
 
 We have the following immediate corollary.
 \begin{corollary}
 Suppose ${\sf prov}^N_\alpha$ is a Fefermanian provability predicate.
 Then, ${\sf prov}^N_\alpha$ is an LR-predicate for $T$ iff  there is a Fefermanian predicate ${\sf prov}_\gamma$ such that
 $T \vdash \opr^N_\gamma B \to \opr^N_\alpha B$, for all $T$-sentences $B$, and ${\sf prov}^N_\gamma$ satisfies L\"ob's Logic
 over $T$.
 \end{corollary}

 We also have:
 
 \begin{corollary}
 Suppose $\opr_\alpha$ is an LR-predicate for $T$.
 Then $T$ interprets $T+\opr_\alpha\bot$. \end{corollary}
 
 \begin{proof}
  Since we have the Interpretation Existence Lemma for Fefermanian provability predicates, we have the result for
 ${\sf prov}_{\alpha^\grullet}$. Moreover $\opr_{\alpha^\grullet}\bot$ is equivalent over $T$ to
 $\opr^\grullet_\alpha\bot$, which is again equivalent to $\opr_\alpha \bot$. 
 By the usual argument, we have $T$ interprets $T+\opr_{\alpha^\grullet}\bot$.
 Hence, $T$ interprets $T+\opr_{\alpha}\bot$.
 \end{proof}
 
 \subsection{Craigification and Smoothening}\label{smoothsmurf}
 To prepare the reader for the example of Section~\ref{aritexa}, we briefly discuss Craigification and
 smoothening. The main point is that smoothening does preserve the Feferman property of being a
 provability predicate of an axiom class.
 
 In this subsection, we follow the sloppy ways of the proof of Theorem~\ref{preserveersmurf}
 and use meta-notations in the object language. We also will suppress the superscript $N$
 that signals relativisation to the chosen numbers.
 
 Consider $\alpha(x)$ of the form $\exists y\, \alpha_0(y,x)$.
 We can transform  $\alpha$ to its Craigification $\alpha^{\sf cr}$ as follows:
 \begin{itemize}
 \item
 $\alpha^{\sf cr}(A)$ iff, for some $n<A$ and $B<A$, we have $A = (\overline n=\overline n \wedge B)$ and
 $\alpha_0(n,B)$.
 \end{itemize}
 
 The smoothening ${\sf prov}_\alpha^{\,\sf sm}$ is defined as follows:
 \begin{itemize}
 \item
 ${\sf prov}_\alpha^{\,\sf sm}(A)$ iff $\exists p\, \exists f\, ({\sf proof}_0(p, A) \wedge \forall B\in {\sf ass}(p)\, \alpha_0(f(B),B))$. \\
 Here `$f$' ranges over finite functions coded as numbers.
 \end{itemize}
 We note that the smoothening takes the syntactic form of $\alpha$ as input, so the notation is a bit misleading.
 The basic insight on the relationship between Craigification and smoothening is simply this:

\begin{theorem}
$T \vdash \forall y\, ({\sf prov}_\alpha^{\,\sf sm}(y) \iff {\sf prov}_{\alpha^{\sf cr}}(y))$.
\end{theorem}

In other words, smoothening preserves the Fefermanian character of a provability predicate.

We will not give the proof here. The main thing is seeing that given an $\alpha$-proof $p$ and the finite function $f$,
we can construct in p-time the corresponding $\alpha^{\,\sf cr}$-proof $p^\ast$.
Conversely, from $p^\ast$ we can efficiently find both a $p$ and and $f$.\footnote{We choose our coding
in such a way that these transformations are feasible. Note that if we had e.g. the Ackermann coding for finite
functions and a string style coding for sequences all this would get far less clear.}

\subsection{An Example}\label{aritexa}
The reader does have to glance through Section 6.2 of
\cite{viss:look19} to understand what is going on in this subsection. 

We use $\beta$ as a standard representation of some single axiom that axiomatises {\sf EA}.
In the example we constructed an axiomatisation $\sigma$ of {\sf EA} with various properties. 
We will show that  ${\sf prov}_\beta$ and ${\sf prov}_\sigma^{\,\sf sm}$ and ${\sf prov}_\sigma^{\grullet}$ and ${\sf prov}_\sigma$
are pairwise distinct  over {\sf EA}. Distinctness means that {\sf EA} does not prove sameness.

We will also show that the G\"odel sentence of the provability predicate
defined is the example has an explicit representation.

\subsubsection{Inclusions}
We first note that since {\sf EA} is finitely axiomatisable, ${\sf prov}_\sigma$ is a uniform provability predicate for {\sf EA}
and, hence, uniform. So, we do have L\"ob's Rule and hence all the desirable properties of $\opr_\sigma^\grullet$.
We have:
\begin{enumerate}[a.]
\item
${\sf EA} \vdash \forall x\, ({\sf prov}_\beta(x) \to {\sf prov}_{\sigma}^{\,\sf sm}(x))$
\item
${\sf EA} \vdash \forall x\, ({\sf prov}^{\,\sf sm}_\sigma(x) \to {\sf prov}_{\sigma}^\grullet(x))$
\item
${\sf EA} \vdash \forall x\, ({\sf prov}^\grullet_\sigma(x) \to {\sf prov}_{\sigma}(x))$
\end{enumerate}

\begin{proof}
We prove (b). We note that ${\sf prov}^{\,\sf sm}_\sigma(x)$ is $\Sigma^0_1$ and that hence we have
\[{\sf EA} \vdash \forall x\, ({\sf prov}^{\,\sf sm}_\sigma(x) \to
 {\sf prov}_\sigma(\gnum{{\sf prov}^{\,\sf sm}_\sigma(\dot x)} )).\]
Let $A :  \forall x\, ({\sf prov}^{\,\sf sm}_\sigma(x) \to {\sf prov}_{\sigma}^\grullet(x))$. We have:
\begin{eqnarray*}
{\sf EA} + \opr_\sigma A \;\; \vdash \;\;  {\sf prov}^{\,\sf sm}_\sigma(a) & \to &  
({\sf prov}_{\sigma}(a) \wedge  {\sf prov}_\sigma(\gnum{{\sf prov}^{\,\sf sm}_\sigma(\dot x)} )) \\
& \to & 
({\sf prov}_{\sigma}(a) \wedge  {\sf prov}_\sigma(\gnum{{\sf prov}^\grullet_\sigma(\dot x)} )) \\
& \to & {\sf prov}^\grullet_{\sigma}(a)
\end{eqnarray*}
It follows that ${\sf EA} \vdash \opr_\sigma A \to A$ and, hence, by L\"ob's rule, ${\sf EA} \vdash A$.
\end{proof}

\subsubsection{Separations}

\emph{We separate ${\sf prov}_\beta$ from ${\sf prov}_\sigma^{\sf sm}$.}
Clearly, ${\sf EA}+\mathcal S^\star\vdash \opr_\sigma^{\,\sf sm} \ell_{\mathfrak p} \neq \overline 1$ and
${\sf EA}+\mathcal S^\star\vdash \opr_\beta \ell_{\mathfrak p} \neq \overline 1 \to \opr_\beta\bot$.
However, as is shown in \cite{viss:look19}, ${\sf EA}+\mathcal S^\star\nvdash \opr_\beta\bot$.
So, ${\sf EA}+\mathcal S^\star\nvdash \opr_\beta\ell_{\mathfrak p} \neq \overline 1$.

\emph{We separate ${\sf prov}_\sigma^{\sf sm}$ from ${\sf prov}_\sigma^\grullet$.}
In ${\sf EA}+\mathcal S^\star$,
smoothening does allow us to use more and more of the
non-standardly finitely many axioms but never all in one proof. It now follows by a minor adaptation of the argument
for Lemma 6.13(b) of \cite{viss:look19} that, over ${\sf EA}+\mathcal S^\star$, we have that
 $\opr_\sigma^{\,\sf sm}\bot$ is equivalent to $\opr_\beta\bot$.
At the same time, $\opr_\sigma^\grullet\bot$ is equivalent to $\opr_\sigma\bot$. Moreover,
by the results of  \cite{viss:look19}, we have the equivalence of $\opr_\sigma \bot$ and $\opr_\beta\opr_\beta\bot$. 
By a model theoretic argument analogous to the argument in \cite{viss:look19}, we find that
${\sf EA}+\mathcal S^\star\nvdash \opr_\beta\opr_\beta \bot \to \opr_\beta \bot$.
So, ${\sf EA}+\mathcal S^\star$ does not prove the equivalence of $\opr_\sigma^{\,\sf sm}\bot$ and $\opr_\sigma^\grullet\bot$.

\emph{We separate ${\sf prov}_\sigma^\grullet$ from ${\sf prov}_\sigma$.}
Suppose we would have the equivalence of 
 ${\sf prov}_\sigma$ and ${\sf prov}_\sigma^{\grullet}$ over  {\sf EA}. 
It would follow that we have L\"ob's Logic for $\opr_\sigma$. But this was refuted in
 \cite{viss:look19}.

\subsubsection{Explicit G\"odel Sentence}
We have seen that the G\"odel sentence for $\opr_\sigma$ is unique over
{\sf EA}. But what is it? The next theorem answers this question.

\begin{theorem}
The G\"odel sentence of $\opr_\sigma$ is, modulo {\sf EA}-provable equivalence, $\oco_\sigma\opr_\sigma\bot$.
\end{theorem}

\begin{proof}
Since we already have uniqueness, it is sufficient to show  that
$\oco_\sigma\opr_\sigma\bot$ is {\sf EA}-provably equivalent to
$\neg\,\opr_\sigma\oco_\sigma\opr_\sigma\bot$. In other words, we want to show that
$\opr_\sigma\oco_\sigma\top$ is {\sf EA}-provably equivalent to $\opr_\sigma\oco_\sigma\opr_\sigma\bot$. 

We remind the reader of the following lemmas of \cite{viss:look19}.
\begin{description}
\item[Lemma 6.9]
${\sf EA}+\mathcal S^\star \vdash \opr_\beta\neg\, \mathcal S^\star$.
\item[Lemma 6.11]
${\sf EA}+\mathcal S^\star \vdash \opr_\sigma A \iff \opr_\beta(\oco_\beta\top \to A)$.
\item[Lemma 6.14(a)] 
${\sf EA}+\mathcal S^\star \vdash \opr_\sigma\oco_\sigma\top$.
\item[Lemma 6.14(b)] 
${\sf EA}+\neg\,\mathcal S^\star \vdash \opr_\sigma\oco_\sigma \top \iff \opr_\beta\bot$.
\end{description}

By Lemma 6.14(b),
it follows that: 
\begin{eqnarray*}
{\sf EA}+\neg\,\mathcal S^\star \vdash\; \opr_\sigma\oco_\sigma\opr_\sigma\bot & \to & \opr_\sigma\oco_\sigma\top \\
& \to & \opr_\beta \bot \\
& \to & \opr_\sigma\oco_\sigma\opr_\sigma\bot
\end{eqnarray*} 
So, in ${\sf EA}+\neg\,\mathcal S^\star$ we have the desired equivalence.

We show we also have the equivalence in ${\sf EA}+\mathcal S^\star$.
In the light of Lemma 6.14(a), it suffices to show  ${\sf EA}+\mathcal S^\star \vdash \opr_\sigma\oco_\sigma\opr_\sigma\bot$.
We reason in ${\sf EA}+\mathcal S^\star$.
By 6.11,  $\opr_\sigma\oco_\sigma\opr_\sigma\bot$ is equivalent to
$ \opr_\beta(\oco_\beta\top\to \oco_\sigma\opr_\sigma\bot)$. By 6.9, we have (i) $\opr_\beta\neg\, \mathcal S^\star$.
By applying necessitation to 6.14(b), we have \[\text{(ii)}\;\; \opr_\beta(\neg\, \mathcal S^\star \to (\opr_\sigma\oco_\sigma \top \iff \opr_\beta\bot)).\]
By combining (i) and (ii), we find that
$\opr_\sigma\oco_\sigma\opr_\sigma\bot$ is equivalent to $ \opr_\beta(\oco_\beta\top\to \oco_\beta\top)$, which in its turn
is equivalent to $\top$. So, we are done.
\end{proof}

\begin{question}
Can we give an example of a modalised fixed point that has no definable solution for the case of ${\sf EA},\sigma$?

Can we given an example of a pair ${\sf EA},\tau$, where $\tau$ is $\Sigma^0_1$ and $\opr_\tau$ is a provability predicate for
{\sf EA} such that the G\"odel sentence is not explicitly definable? If not, what about the more general case?
\end{question}

\begin{question}
Consider an $\Sigma^0_1$-predicate $\tau$ that axiomatises {\sf EA} in {\sf EA}. The provability logic of $\opr_\tau^\grullet$ contains
{\sf GL}. However, \emph{prima facie}, Solovay's proof fails. Can we still prove that the logic is precisely {\sf GL}.
What about the logic for $\opr_\sigma^\grullet$ for the specific predicate $\sigma$ studied above?
\end{question}

\end{document}